\documentclass[11pt]{article}
\usepackage[utf8]{inputenc}
\usepackage{amsmath,amssymb,bbm}
\usepackage{stmaryrd}
\usepackage{comment}
\usepackage[dvips]{color}
\usepackage[T1]{fontenc}

\usepackage{pdfpages,graphics}
\usepackage{todonotes}
\usepackage{enumitem}
\setlist{nolistsep}
\setlist[itemize]{labelindent=*,leftmargin=*}
\setlist[enumerate]{labelindent=*,leftmargin=*}

\newtheorem{defi}{Definition}
\newtheorem{prop}[defi]{Proposition}
\newtheorem{theo}[defi]{Theorem}
\newtheorem{conj}[defi]{Conjecture}
\newtheorem{lemm}[defi]{Lemma}
\newtheorem{coro}[defi]{Corollary}
\newtheorem{rema}[defi]{Remark}
\newtheorem{exem}[defi]{Example}
\newtheorem{exems}[defi]{Examples}

\newcommand{\bdefi}{\begin{defi}}
\newcommand{\edefi}{\end{defi}}
\newcommand{\bprop}{\begin{prop}}
\newcommand{\eprop}{\end{prop}}
\newcommand{\btheo}{\begin{theo}}
\newcommand{\etheo}{\end{theo}}
\newcommand{\blemm}{\begin{lemm}}
\newcommand{\brema}{\begin{rema}}
\newcommand{\erema}{\end{rema}}
\newcommand{\bexer}{\begin{exem}}
\newcommand{\eexer}{\end{exem}}
\newcommand{\bexems}{\begin{exems}}
\newcommand{\eexems}{\end{exems}}
\newcommand{\bconj}{\begin{conj}}
\newcommand{\econj}{\end{conj}}
\newcommand{\elemm}{\end{lemm}}
\newcommand{\bcoro}{\begin{coro}}
\newcommand{\ecoro}{\end{coro}}

\usepackage{mathrsfs}
\renewcommand\mathcal{\mathscr}

\newcommand{\N}{{\cal N}}

\newcommand{\F}{{\cal F}}
\newcommand{\V}{{\cal V}}
\newcommand{\W}{{\cal W}}

\renewcommand{\H}{{\cal H}}
\newcommand{\OOO}{{\cal O}}
\newcommand{\C}{{\cal C}}

\newcommand{\Scal}{{\cal S}}
\newcommand{\maths}[1]{{\mathbb #1}}  

\newcommand{\RR}{\maths{R}}
\newcommand{\NN}{\maths{N}}
\newcommand{\CC}{\maths{C}}
\newcommand{\QQ}{\maths{Q}}

\newcommand{\SSS}{\maths{S}}

\newcommand{\HH}{\maths{H}}

\newcommand{\ZZ}{\maths{Z}}
\newcommand{\PP}{\maths{P}}

\newcommand{\aaa}{{\mathfrak a}}

\newcommand{\weakstar}{\overset{*}\rightharpoonup}
\newcommand{\ra}{\rightarrow}
\newcommand{\bs}{\backslash}

\newcommand{\wt}[1]{{\widetilde{#1}}}

\newcommand{\ga}{\gamma}
\newcommand{\Ga}{\Gamma}

\newcommand{\card}{{\operatorname{Card}}}
\renewcommand{\Re}{{\operatorname{Re}}}
\renewcommand{\Im}{{\operatorname{Im}}}
\newcommand{\Vol}{\operatorname{Vol}}
\newcommand{\Isom}{\operatorname{Isom}}

\newcommand{\covol}{\operatorname{Covol}}

\newcommand{\PSL}{\operatorname{PSL}}
\newcommand{\SL}{\operatorname{SL}}

\newcommand{\bigO}{\operatorname{O}}

\newcommand{\Heis}{\operatorname{Heis}}

\newcommand{\haarheis}{\operatorname{Haar}_{\operatorname{Heis}_3}}
\renewcommand{\log}{\operatorname{ln}}
\newcommand{\Leb}{\operatorname{Leb}}

\newcommand{\Sp}{\operatorname{Sp}}
\newcommand{\WA}{\operatorname{WA}}

\newcommand{\PSLZ}{\operatorname{PSL}_{2}(\ZZ)}

\newcommand{\PSU}{\operatorname{PSU}}

\newcommand{\tr}{\operatorname{\tt tr}}

\newcommand{\diam}{\operatorname{diam}}

\newcommand{\flow}[1]{{\phi_{#1}}}

\newcommand\Perp{\operatorname{Perp}}

\newcommand\normalout{\partial^1_{+}}
\newcommand\normalin{\partial^1_{-}}
\newcommand\normalpm{\partial^1_{\pm}}
\newcommand\normalmp{\partial^1_{\mp}}

\newcounter{fig}

\title{A survey of some arithmetic applications 
\\ of ergodic theory in negative curvature}
\author{Jouni Parkkonen \and Fr\'ed\'eric Paulin} 
\date{\today}

\begin{document}
\bibliographystyle{../alphanum}
\maketitle

\begin{abstract}
  This paper is a survey of some arithmetic applications of techniques
  in the geometry and ergodic theory of negatively curved Riemannian
  manifolds, focusing on the joint works of the authors. We describe
  Diophantine approximation results of real numbers by quadratic
  irrational ones, and we discuss various results on the
  equidistribution in $\RR$, $\CC$ and in the Heisenberg groups of
  arithmetically defined points. We explain how these results are
  consequences of equidistribution and counting properties of common
  perpendiculars between locally convex subsets in negatively curved
  orbifolds, proven using dynamical and ergodic properties of their
  geodesic flows. This exposition is based on lectures %notes
  at the conference ``Chaire Jean Morlet: Géométrie et systèmes
  dynamiques'', at the CIRM, Luminy, 2014.  We thank B.~Hasselblatt
  for his strong encouragements to write this survey. 
\end{abstract}

\section{Introduction}
\label{sect:intro}

For several decades, tools from dynamical systems, and in particular
ergodic theory, have been used to derive arithmetic and number
theoretic, in particular Diophantine approximation results, see for
instance the works of Furstenberg, Margulis, Sullivan, Dani,
Kleinbock, Clozel, Oh, Ullmo, Lindenstrauss, Einsiedler, Michel,
Venkatesh, Marklof, Green-Tao, Elkies-McMullen, Ratner, Mozes, Shah,
Gorodnik, Ghosh, Weiss, Hersonsky-Paulin, Parkkonen-Paulin and many
others, and the references \cite{Kleinbock01, Lindenstrauss07,
  Kleinbock09, ArbMatMor09, Athreya09, GorNev10, EinWar11, ParPauRev}.

In Subsection \ref{subsec:approxframwork} of this survey, we introduce
a general framework of Diophantine approximation in measured metric
spaces, in which most of our arithmetic corollaries are inserted (see
the end of Subsection \ref{subsec:approxframwork} for references
concerning this framework). In order to motivate it, we first recall
in Subsection \ref{subsec:bsaciDiophApprox} some very basic and
classical results in Diophantine approximation (see for instance
\cite{Bugeaud04,Bugeaud12}). A selection (extracted from
\cite{ParPau11MZ,ParPau14AFST,ParPauHeis}) of our arithmetic results are
then stated in Subsections \ref{subsec:DiophapproxR},
\ref{subsec:equidistribRC} and \ref{subsec:equidistribHeisenberg},
where we indicate how they fit into this framework: Diophantine
approximation results (\`a la Khintchine, Hurwitz, Cusick-Flahive and
Farey) of real numbers by quadratic irrational ones, equidistribution
of rational points in $\RR$ (for various height functions), in $\CC$
and in the Heisenberg group, ...

We will explain in Subsection \ref{subsect:fromgeomtoarith} the
starting point of their proofs, using the geometric and ergodic tools
and results previously described in Subsection
\ref{subsect:commonperp}, where we give an exposition of our work in
\cite{ParPau13b}: an asymptotic formula as $t\ra+\infty$ for the
number of common perpendiculars of length at most $t$ between closed
locally convex subsets $D^-$ and $D^+$ in a negatively curved
Riemannian orbifold, and an equidistribution result of the initial and
terminal tangent vectors $v^-_\alpha$ and $v^+_\alpha$ of the common
perpendiculars $\alpha$ in the outer and inner unit normal bundles of
$D^-$ and $D^+$, respectively.  Common perpendiculars have been
studied, in various particular cases, sometimes not explicitly, by
Basmajian, Bridgeman, Bridgeman-Kahn, Eskin-McMullen, Herrmann, Huber,
Kontorovich-Oh, Margulis, Martin-McKee-Wambach, Meyerhoff, Mirzakhani,
Oh-Shah, Pollicott, Roblin, Shah, the authors and many others (see the
comments after Theorem \ref{theo:maincount} below, and the survey
\cite{ParPauRev} for references). 

Section \ref{sect:measures} presents the background notions on the
geometry in negative curvature, describes various useful measures, and
recalls the basic results about them, due to works of Patterson,
Sullivan, Bowen, Margulis, Babillot, Roblin, Otal-Peigné,
Kleinbock-Margulis, Clozel, Oh-Shah, Mohammadi-Oh and the authors (see
for instance \cite{Roblin03, Babillot02b, OtaPei04, MohOh14,
  ParPau14ETDS}).  See \cite{PauPolSha,BroParPau15} for extensions to
manifolds with potentials and to trees with potentials.

\medskip
Let us denote by $^cA$ the complementary subset of a subset $A$ of any
given set, by $\|\mu\|$ the total mass of a measure $\mu$, by
$\Leb_\RR$ and $\Leb_\CC$ the Lebesgue measures on $\RR$ and $\CC$, by
$\Delta_x$ the unit Dirac mass at any point $x$ in any topological
space, and by $\stackrel{*}{\rightharpoonup}$ the weak-star
convergence of measures on any locally compact space.

\section{Arithmetic applications}
\label{sect:arithmeticapplications}

\subsection{Basic and classic Diophantine approximation }
\label{subsec:bsaciDiophApprox}

When denoting a rational number $\frac{p}{q}\in\QQ$, we will assume
that $p$ and $q$ are coprime and that $q>0$. For every irrational real
number $x\in\RR-\QQ$, let us define the {\it approximation exponent}
$\omega(x)$ of $x$ as
$$
\omega(x)=\limsup_{\frac{p}{q}\in\QQ,\;
  q\ra+\infty}\frac{-\ln\big|x-\frac{p}{q}\big|}{\ln q}\;.
$$
The Dirichlet theorem implies that 
$$
\inf_{x\in\RR-\QQ}w(x)=2\,, 
$$ 
which motivates the definition of the  
the {\it approximation constant} $c(x)$ of $x\in\RR-\QQ$ as
$$
c(x)=\liminf_{\frac{p}{q}\in\QQ,\; q\ra+\infty} q^2\Big|x-\frac{p}{q}\Big|\;.
$$
The convention varies, some other references consider $c(x)^{-1}$ or
$(2c(x))^{-1}$ as the approximation constant.  The {\it Lagrange
  spectrum} for the approximation of real numbers by rational ones is
$$
\Sp_{\rm Lag}=\{c(x)\;:\; x\in\RR-\QQ\}\;.
$$
Given $\psi:\NN\ra \mathopen{]}0,+\infty\mathclose{[}$, the set of
{\it $\psi$-well approximable} real numbers by rational ones is$$
\WA_\psi=\big\{x\in\RR-\QQ\;:\; 
\card\big\{\frac{p}{q}\in\QQ\;:\;\Big|x-\frac{p}{q}\Big|\leq 
\psi(q)\big\}=+\infty\big\}\;.
$$
Again the convention varies, some other references consider $q\mapsto
\psi(q)/q$ or similar instead of $\psi$.

Denoting by $\dim_{\rm Hau}$ the Hausdorff dimension of subsets of
$\RR$, the Jarnik-Besicovich theorem says that, for every $c\geq 0$,
$$
\dim_{\rm Hau}\{x\in\RR-\QQ\;:\; w(x)\geq 2+c\}=\frac{2}{2+c}
$$
Many properties of the Lagrange spectrum are known (see for instance
\cite{CusFla89}): $\Sp_{\rm Lag}$ is bounded, with maximum
$\frac{1}{\sqrt{5}}$ known as the {\it Hurwitz constant}
(Korkine-Zolotareff 1873, Hurwitz 1891); it is closed (Cusick 1975),
and contains a maximal interval $[0,\mu]$ with $\mu>0$ (Hall 1947),
with in fact $\mu=491993569/(2221564096+283748\sqrt{462})\simeq
0.2208$ known as the {\it Freiman constant} (Freiman 1975).

The Khintchine theorem gives a necessary and sufficient criteria for
the $\psi$-well approximable real number to have full or zero Lebesgue
measure:
$$
\left\{\begin{array}{ll}
\Leb_\RR({}^c\WA_\psi)=0& 
{\rm if~} \sum_{q=1}^{+\infty} \,q\,\psi(q)=+\infty
\\ \Leb_\RR(\WA_\psi)=0& 
{\rm if~} \sum_{q=1}^{+\infty} \,q\,\psi(q)<+\infty\;.
\end{array}\right.
$$
Finally, the equidistribution of Farey fractions (which is closely
related with the Mertens formula) is an equidistribution theorem of
the rational numbers in $\RR$ when their denominator tends to
$+\infty$:
$$\frac{\pi^2}{12\,s}\sum_{\frac{p}{q}\in\QQ,\;|q|\leq s} 
\Delta_{\frac{p}{q}}\;\; \stackrel{*}{\rightharpoonup}\;\;\Leb_\RR\;.
$$

\subsection{An approximation framework}
\label{subsec:approxframwork}

The general framework announced in the introduction concerns the
quantitative answers to how dense a given dense subset of a given
topological set is.

Let $(Y,d,\mu)$ be a metric measured space (see for instance
\cite[Chap.~3$\frac12$]{Gromov99a} and \cite{Heinonen01} for
generalities), with $Y$ a subspace of a topological space $X$, let $Z$
be a countable (to simplify the setting in this survey) subset of $X$
whose closure contains $Y$ (for instance a dense orbit of a countable
group of homeomorphisms of $X$), and let $H:Z\ra \mathopen{]}0,
+\infty[$ be a map called a {\it height function} which is proper (for
every $r>0$, the set $H^{-1}(\mathopen{]}0,r])$ is finite).  The
classical Diophantine approximation problems of subsection
\ref{subsec:bsaciDiophApprox} fit into this framework with $X=\RR$,
$Y=\RR-\QQ$, $Z=\QQ$ and $H:\frac{p}{q}\mapsto q$ or $H:\frac{p}{q}
\mapsto q^2$, considered modulo translation by integers.  The height
functions in question are invariant under translation by $\ZZ$, and
the properness condition is satisfied in the quotient space
$\ZZ\bs\RR$ or, equivalently, by restriction to the unit interval.

We endow $Z$ with the Fr\'echet filter of the complementary subsets of
its finite subsets: $z\in Z$ tends to infinity if and only if $z$
leaves every finite subset of $Z$. Generalising the definitions of
Subsection \ref{subsec:bsaciDiophApprox}, for every $y\in Y$, we may
define the {\it approximation exponent} $\omega(y)=
\omega_{X,\,Y,\,Z,\,H}(y)$ of $y$ by the elements of $Z$ with height
function $H$ as
$$
\omega(y)=\limsup_{z\in Z}\;\frac{-\ln d(y,z)}{\ln H(z)}\;,
$$
and the {\it approximation constant} $c(y)=c_{X,\,Y,\,Z,\,H}(y)$ of
$y$ by the elements of $Z$ with height function $H$ as
$$
c(y)=\liminf_{z\in Z} \;H(z)\,d(y,z)\;.
$$
The {\it Lagrange spectrum} $\Sp_{\rm Lag}=\Sp_{X,\,Y,\,Z,\,H}$ for the
approximation of the elements of $Y$ by the elements of $Z$ with
height function $H$ is
$$
\Sp_{\rm Lag}=\{c(y)\;:\; y\in Y\}\;.
$$
Its least upper bound will be called the {\it Hurwitz constant} for
the approximation of the elements of $Y$ by the elements of $Z$ with
height function $H$. Given $\psi:\mathopen{]}0,+\infty\mathclose{[}
\ra \mathopen{]}0, +\infty\mathclose{[}$, the set $\WA_\psi=
\WA_{\psi,\,X,\,Y,\,Z,\,H}$ of {\it $\psi$-well approximable} elements
of $Y$ by the elements of $Z$ with height function $H$ is
$$
\WA_\psi=\big\{y\in Y\;:\; 
\card\big\{z\in Z\;:\;d(y,z)\leq 
\psi(H(z))\big\}=+\infty\big\}\;.
$$

The approximation problems may be subdivided into several classes, as
follows, possibly by taking the appropriate height function, for
appropriate maps $\psi:\mathopen{]}0,+\infty\mathclose{[}\ra
\mathopen{]}0,+\infty\mathclose{[}$.
\begin{itemize}
\item The approximation exponent problem: study the map
  $y\mapsto \omega(y)$.
\item The Lagrange problem: study the Lagrange spectrum for
  the approximation of the elements of $Y$ by the elements of $Z$ with
  height function $H$.
\item The Jarnik-Besicovich problem: compute the Hausdorff
  dimension of the set of $y\in Y$ with $c(y)\geq \psi(H(y))$.
\item The Khintchine problem: study whether $\mu$-almost
  every (or $\mu$-almost no) $y\in Y$ is $\psi$-well approximable.
\item The counting problem: study the asymptotics as $s$
  tends to $+\infty$ of 
$$
\card\{y\in Y\;:\; H(y)\leq s\}\;.
$$
\item The equidistribution problem: study the set of
  weak-star accumulation points as $s$ tends to $+\infty$ of the
  probability measures
$$
\frac{1}{\card\{y\in Y\;:\; H(y)\leq s\}}\;\sum_{z\in Z,\;H(z)\leq s}\;\Delta_z\;.
$$
\end{itemize}
The equidistribution problem, closely linked to the counting problem,
is the one we will concentrate on in Subsections
\ref{subsec:DiophapproxR}, \ref{subsec:equidistribRC} and
\ref{subsec:equidistribHeisenberg}.

\medskip This framework is not new (see for instance the works of
Kleinbock), and many results have developped some aspect of it.

\medskip
(1) For instance, $X$ could be the boundary at infinity of a Gromov
hyperbolic metric space, $Y$ a (subset of) the limit set of a discrete
group $\Ga$ of isometries of this hyperbolic space, $Z$ could be the
orbit under $\Ga$ of some point $x\in X$, and $H:Z\ra
\mathopen{]}0,+\infty\mathclose{[}$ could be $\ga x\mapsto 1 +
d_X(A_0,\ga B_0)$ where $A_0,B_0$ are subsets of $X$, with $B_0$
invariant under the stabilizer of $x$ in $\Ga$. Numerous aspects of
this particular case have been developped, by Patterson, Sullivan,
Dani, Hill, Stratmann, Velani, Bishop-Jones, Hersonsky-Paulin,
Parkkonen-Paulin, and the most complete and general version is due to
Fishman-Simmons-Urba\'nski \cite{FisSimUrb14}, to which we refer in
particular for their thorough long list of references.

\medskip (2) The case when $X=\RR^N$, $Y$ is a submanifold (or a more
general subset) of $X$, $Z=\QQ^N$ has been widely studied, under the
name of Diophantine approximation on curves, submanifolds and fractals
subsets, by many authors, including Kleinbock-Margulis, Bernik,
Dodson, Beresnevich, Velani, Kleinbock-Weiss and others, see for
instance \cite{BerDod99, BadBerVel13} and their references.

\medskip (3) If $X=\underline{X}(\RR)$ is for instance the set of real
points of an algebraic manifold defined over $\QQ$, if $Z=
\underline{X}(\QQ)$ is the set of rational points, if $Y$ is the
(Hausdorff) closure of $Z$ in $X$ or this closure minus $Z$, there are
many results, in particular when $X$ is homogeneous, on the above
Diophantine approximation problems. Similarly, if $X=G/H$ is a
homogeneous space of a semisimple connected Lie group $G$, if $Z$ is
the orbit in $X$ of a lattice in $G$ and if $Y$ is the closure of $Z$
in $X$, most of the above problems have been stated and studied for
instance in \cite{GhoGorNev12, GhoGorNev14a, GhoGorNev14b}. We refer
to the works of Benoist, Browning, Colliot-Thélène, Duke, Einsiedler,
Eskin, Gorodnik, Heath-Brown, Lindenstrauss, Margulis, Mozes, Oh,
Ratner, Quint, Rudnick, Salberger, Sarnak, Shah, Tomanov, Ullmo,
Venkatesh, Weiss and many others for equidistribution results in
homogeneous spaces, see for instance \cite{Breuillard05, BjoFis09,
  Harcos10, EinWar11, Serre12, GreTao12, Kim13, BenQui12, BenQui13a}.

\medskip In the next three subsections, we give a sample of the
results obtained by the authors on the above problems.

\subsection{Diophantine approximation in $\RR$ by quadratic
  irrationals}
\label{subsec:DiophapproxR}

Our first results concern the Diophantine approximation  of real
numbers, where we replace the approximating rationals by quadratic
irrational numbers. We refer for instance to \cite{Bugeaud04} and its
references for very different Diophantine approximation results by
algebraic numbers.

We denote by $\alpha^\sigma$ the Galois conjugate of a real quadratic
irrational $\alpha$, and by $\tr \alpha =\alpha+\alpha^\sigma$ its
trace. We approximate the real points by the elements of the orbit of
a fixed quadratic irrational $\alpha_0$ by homographies under
$\PSL_2(\ZZ)$ (and their Galois conjugates). We denote by $\ga\cdot x$
the action by homography of $\ga\in\PSL_2(\RR)$ on $x\in\PP_1(\RR)=
\RR\cup\{\infty\}$.

Let $X=\RR$ (with the standard Euclidean distance and Lebesgue
measure), let $Y=\RR-\QQ-\PSL_2(\ZZ)\cdot\alpha_0$, let $Z=
\PSL_2(\ZZ)\cdot\alpha_0$ and let $H:\alpha\mapsto
\frac{2}{|\alpha-\alpha^\sigma|}$. It turns out that $H$ is an
appropriate height function on each $\PSL_2(\ZZ)$-orbit under
homographies (working modulo translation by $\ZZ$) of a given
quadratic irrational. We refer to \cite{ParPau11MZ,ParPau12JMD} for a
proof of this, and for more algebraic expressions of this height
function and its important differences with classical ones.

We consider in the first statement below the particular case when
$\alpha_0$ is the Golden Ratio $\phi= \frac{1+\sqrt{5}}{2}$, and we
refer to \cite{ParPau11MZ} for the general version. One way this
restriction simplifies the statement is that the Golden Ratio is in
the same orbit under $\PSL_2(\ZZ)$ as its Galois conjugate, which is
not the case of every quadratic irrational (see for instance
\cite{Sarnak07}).  The next result is proven in \cite[Theo.~1.3,
Prop.~1.4]{ParPau11MZ}, with a mistake in the Hurwitz constant
corrected in the erratum of loc. cit., thanks to Bugeaud, who gave
another way to compute it in \cite{Bugeaud14}, using only continued
fractions techniques.

\btheo[Parkkonen-Paulin] \label{theo:Lagrangequadirrat} With the
notation $X,Y,Z,H$ as above, the Lagrange spectrum $\Sp_{\rm Lag}$ is
closed, bounded, with Hurwitz constant $\frac{3}{\sqrt{5}}-1$. For
every $\psi:\mathopen{]}0,+\infty\mathclose{[} \ra \mathopen{]}0,
+\infty\mathclose{[}$ such that $t\mapsto \ln\psi(e^t)$ is Lipschitz,
we have
$$
\left\{\begin{array}{ll} \Leb_\RR({}^c\WA_\psi)=0& {\rm if~}
    \int_1^{+\infty} \,\psi(t)\,dt=+\infty \\ \Leb_\RR(\WA_\psi)=0&
    {\rm if~} \int_1^{+\infty} \,\psi(t)\,dt<+\infty\;.
\end{array}\right.
$$
\etheo

The exact value of the Hurwitz constant for the Diophantine
approximation of real numbers by elements of an orbit under
$\PSL_2(\ZZ)$ (or a congruence subgroup) of a general quadratic
irrational (and its Galois conjugate) is an interesting open problem.

\medskip 
%For every quadratic irrational $\alpha$, we denote by $\tr
%\alpha =\alpha+\alpha^\sigma$ its trace, by $\n(\alpha)= \alpha
%\alpha^\sigma$ its norm, and by $z\alpha:t\mapsto t^2+(\tr \alpha)t +
%\n(\alpha)$ the monic rational quadratic polynomial with roots $\alpha$
%and $\alpha^\sigma$.
%
Let $\alpha_0,\beta_0$ be fixed integral quadratic irrationals, and
let $R_{\alpha_0},R_{\beta_0}$ be the regulators of the lattices $\ZZ+
\ZZ\alpha_0,\ZZ+ \ZZ\beta_0$ respectively. The integrality assumption
is only present here in order to simplify the statements below in this
survey, see \cite{ParPau14AFST} for the general version.

The following result is an equidistribution result of the traces of
the quadratic irrationals in a given orbit (by homographies) under
$\PSL_2(\ZZ)$ of a quadratic irrational, using the above height
function $H:\alpha\mapsto \frac{2}{|\alpha-\alpha^\sigma|}$. We refer
to \cite[Theo.~4.1]{ParPau14AFST} for a version with additional
congruence assumptions, and to \cite[Theo.~4.2]{ParPau14AFST} and
\cite[Theo.~4.4]{ParPau14AFST} for extensions to quadratic irrationals
over an imaginary quadratic number field (using relative traces) or
over a rational quaternion algebra.

\btheo[Parkkonen-Paulin]\label{theo:equidistrace} As $s\ra +\infty$,
we have
$$
\frac{\pi^2}
{3\,R_{\alpha_0}\,s }\;\;
\sum_{\alpha\in \PSL_2(\ZZ) \cdot \alpha_0
\;:\; H(\alpha)\leq s}\Delta_{\tr\, \alpha}\;\;
\stackrel{*}{\rightharpoonup}\;\; \Leb_\RR\,.
$$
\etheo 

We introduced in \cite{ParPau14AFST} another height function for the
Diophantine approximation of real numbers by the elements of the orbit
(by homographies) under $\PSL_2(\ZZ)$ of $\beta_0$, which measures
their relative complexity with respect to $\alpha_0$. Let 
$$[a,b,c,d]
=\frac{(c-a)\,(d-b)}{(c-b)\,(d-a)}$$ be the standard crossratio of a
quadruple $(a,b,c,d)$ of pairwise distinct points in $\RR$. For every
$\beta\in \PSLZ\cdot\beta_0- \{\alpha_0,\alpha_0^\sigma\}$, let
$$
H_{\alpha_0}(\beta)=
\frac{1}{\max\{|[\alpha_0,\beta,\alpha_0^\sigma,\beta^\sigma]|,\;
  |[\alpha_0,\beta^\sigma,\alpha_0^\sigma,\beta]|\}}\;.
$$
We prove in \cite{ParPau14AFST} that the map $H_{\alpha_0}$ is an
appropriate height function modulo the action (by homographies) of the
fixator $\PSLZ_{\alpha_0}$ of $\alpha_0$ in $\PSLZ$.  The
following theorem is a counting result of quadratic irrationals relative to a given
one.

\btheo[Parkkonen-Paulin]\label{theo:relatcomplcount} There exists
$\kappa>0$ such that, as $s\ra +\infty$,
$$
\card \{\beta\in \PSLZ_{\alpha_0}\bs
\PSLZ\cdot\beta_0,\; H_{\alpha_0}(\beta)\le s\} =
\frac{48\,R_{\alpha_0}\,R_{\beta_0}}
{\pi^2}\; s + \bigO(s^{1-\kappa})\,.
$$
\etheo

We refer to \cite[Theo.~4.9]{ParPau14AFST} for a more general version,
including additional congruence assumptions, and to
\cite[Theo.~4.10]{ParPau14AFST} for an extension to quadratic irrationals
over an imaginary quadratic extension of $\QQ$.  Theorem
\ref{theo:relatcomplcount} fits into the framework of Subsection
\ref{subsec:approxframwork} with $ X=Y=\PSLZ_{\alpha_0}\bs(\PP_1(\RR)-
\{\alpha_0,\alpha_0^\sigma\})$,
$Z=\PSLZ_{\alpha_0}\bs(\PSLZ\cdot\beta_0-
\{\alpha_0,\alpha_0^\sigma\}) $ and $H:\PSLZ_{\alpha_0}\cdot\beta
\mapsto H_{\alpha_0}(\beta)$.

\subsection{Equidistribution of rational points in $\RR$ and $\CC$}
\label{subsec:equidistribRC}

In this section, we will consider equidistribution results in $\RR$
and $\CC$ of arithmetically defined points.

First, we would like to consider again the approximation of points in
$\RR$ by points in $\QQ$, but to change the height function $H$, using
an indefinite rational binary quadratic form $Q$ which is not the
product of two rational linear forms, by taking
$$
H\big(\frac{p}{q}\big)=|Q(p,q)|.
$$ 
It is easy to see that this is (locally) an appropriate height
function outside the roots $\alpha,\alpha^\sigma$ of $t\mapsto
Q(t,1)$: for every $s\geq 0$, the number of $\frac{p}{q}\in\QQ$ such
that $H\big(\frac{p}{q}\big)\leq s$ is locally finite in
$\RR-\{\alpha, \alpha^\sigma\}$. This fits in the framework of
Subsection \ref{subsec:approxframwork} as follows. Let
$$
\operatorname{SO}_Q(\ZZ)=\{g\in \SL_2(\ZZ)\;:\; Q\circ g=Q\}
$$ 
be the integral group of automorphs of $Q$. Note that the subgroup
$\operatorname{SO}_Q(\ZZ)$ of $\SL_2(\ZZ)$ injects into $\PSL_2(\ZZ)$
if and only if $\tr \alpha \neq 0$. Let $\operatorname{PSO}_Q(\ZZ)$ be
the image of $\operatorname{SO}_Q(\ZZ)$ in $\PSL_2(\ZZ)$, which acts
by homographies on $\PP_1(\RR)$, fixing $\alpha$ and $\alpha^\sigma$,
acting properly discontinuously on $\PP_1(\RR)-\{\alpha,
\alpha^\sigma\}$.  We therefore study the Diophantine approximation
problems with $X=\;\operatorname{PSO}_Q(\ZZ) \bs(\PP_1(\RR)- \{\alpha,
\alpha^\sigma\})$ and $Z$ the image of $\PP_1(\QQ)$ by the canonical
projection $(\PP_1(\RR)- \{\alpha, \alpha^\sigma\})\ra X$.

We only consider in this survey the particular case when $Q$ is
$Q(u,v)=u^2-uv-v^2$ (closely related to the Golden Ratio, as the
knowledgeable reader has already seen!). We refer to
\cite[Theo.~5.10]{ParPau14AFST} for the general result, with error term
estimates and additional congruence assumptions, and to
\cite{GorPau15} for an extension to $n$-ary norm forms for general
$n>2$.

\btheo[Parkkonen-Paulin]\label{theo:equidistribR} As $s\ra+\infty$, we
have
$$
\frac{\pi^2}{6\;s}\;
\sum_{\frac{p}{q}\in\QQ,\;|p^2-pq -q^2|\leq s} 
\Delta_{\frac pq}
\;\;\weakstar\;\;\frac{dt}{|t^2-t-1|}\,.
$$
\etheo

Note that the measure to which the rational points equidistribute when
counted with this multiplicity is no longer the Lebesgue measure, but
is the natural smooth measure invariant under the real group of
automorphs of $Q$ (unique up to multiplication by a locally constant
positive function).

\medskip 
Now, let us turn to the Diophantine approximation of complex numbers
by Gaussian rational ones. Every element of the imaginary quadratic
field $K=\QQ(i)$ of Gaussian rational numbers may and will be written
$\frac pq$ with $p,q\in \ZZ[i]$ relatively prime Gaussian integers,
and this writing is unique up to the multiplication of $p$ and $q$ by
the same invertible Gaussian integer $1$, $-1$, $i$ or $-i$. In
particular, the map
$$
H\big(\frac{p}{q}\big)=|q|
$$ 
is well defined, and is clearly an appropriate height function on
$\QQ(i)$ modulo $\ZZ[i]$. We hence consider, with the notation of
Subsection \ref{subsec:approxframwork}, the spaces $X=\CC/\ZZ[i]$,
$Z=\QQ(i)/\ZZ[i]$, $Y=(\CC-\QQ(i))/\ZZ[i]$ and the above height
function $H:Z\ra[0,+\infty[\,$.

The following result (due to \cite[Coro.~6.1]{Cosentino99} albeit in a
less explicit form) is an equidistribution result of the Gaussian
rational points in the complex field, analogous to the Mertens theorem
on the equidistribution of Farey fractions in the real field. We
denote, here and in Subsection \ref{subsec:equidistribHeisenberg}, by
$\OOO_K=\ZZ[i]$ the ring of integers of $K=\QQ(i)$, by
$\OOO_K^\times=\{1, -1, i,-i\}$ the group of invertible elements in
$\OOO_K$, by $D_K=-4$ its discriminant, and by
$$
\zeta_{K}:s \mapsto
\sum_{\aaa {\rm ~nonzero~ideal~in~} \OOO_K} \frac1{N(\aaa)^s}
$$ 
its Dedekind zeta function. The pictures below show the fractions
$\frac pq\in\QQ[i]$ in the square $[-1,1]\times[-1,1]$ where $p,q\in
\ZZ[i]$ are relatively prime Gaussian integers with $|q|\leq
5$ and $|q|\leq
10$. The fact that there is a large white region around the
fractions $\frac pq\in\QQ[i]$ with $|q|$ small will be explained in
Subsection 4.2.

\begin{center}
\includegraphics[width=7cm]{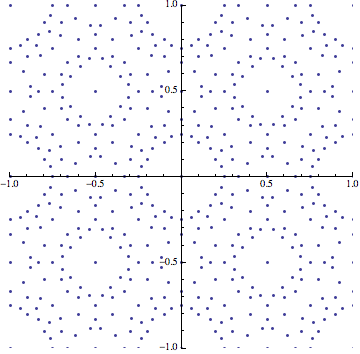}
\includegraphics[width=7cm]{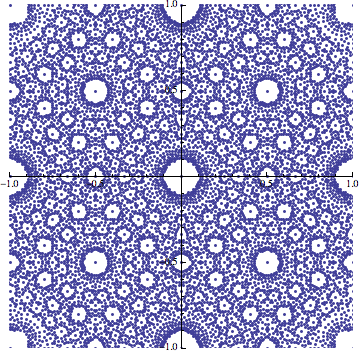}
\end{center}

\btheo[Cosentino, Parkkonen-Paulin]\label{theo:equidistribC} 
As $s\to+\infty$,
$$
\frac{|\OOO_K^\times|\,|D_K|\,\zeta_{K}(2)}{2\pi\,  s^2}
\sum_{\frac pq\in K\;:\;|q|\leq s} \Delta_{\frac pq}\;\;
\weakstar\;\;\Leb_\CC\,.
$$
\etheo

We refer to \cite[Theo.~1.1]{ParPau14AFST} for a version of this theorem
valid for any imaginary quadratic number field (see below a plot in
the square $[-1,1]\times[-1,1]$ of the fractions $\frac pq\in
\QQ[i\sqrt{3}]$ where $p,q\in \ZZ[e^{\frac{2\pi}{3}}]$ are relatively
prime Eisenstein integers with $|q|\leq 5$ and $|q|\leq 10$),  
with additional congruence assumptions on $u$ and $v$, to
\cite[Theo.~1.2]{ParPau14AFST} for analogous results in Hamilton's
quaternion division algebra, and to \cite{Cosentino99, ParPau14AFST} for
error term estimates.

\begin{center}
\includegraphics[width=7cm]{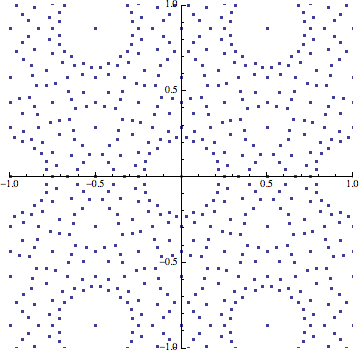}
\includegraphics[width=7cm]{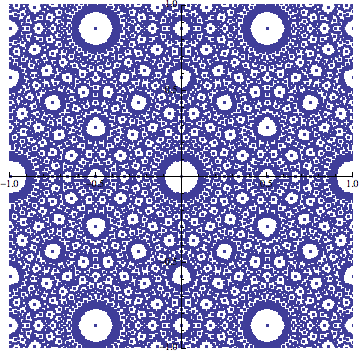}
\end{center}

\subsection{Equidistribution and counting in the 
Heisenberg group}
\label{subsec:equidistribHeisenberg}

In this section, we will consider Diophantine approximation results of
elements of the Heisenberg group by arithmetically defined points.

Recall that the $3$-dimensional {\it Heisenberg group} is the
nilpotent real Lie group with underlying manifold
$$
\Heis_3=\{(w_0,w)\in \CC\times \CC\;:\; 2\Re\;w_0=|w|^2\}
$$
and law $(w_0,w)(w'_0,w')= (w_0+w'_0+w'\,\overline{w},w+w')$. A
standard model in control theory of $\Heis_3$, with underlying
manifold $\RR^3$, may be obtained by the change of variables
$(x,y,t)\in\RR^3$ with $w=x+iy$ and $t=2\,\Im\;w_0$.  We endow
$\Heis_3$ with its Haar measure
$$
d\haarheis(w_0,w)= d(\Im\,w_0)\, d(\Re\,w)\, d(\Im\,w)\,,
$$
(that is $\frac{1}{2}\,dx\,dy\,dt$ in the above coordinates
$(x,y,t)\,$) and with the (almost) distance $d''_{Cyg}$ defined below.

The {\it Cygan distance} $d_{\rm Cyg}$ on $\Heis_3$ (see \cite[page
160]{Goldman99}, sometimes called the {\it Kor\'anyi distance} in
sub-Riemannian geometry, though Kor\'anyi \cite{Koranyi85} attributes
it to Cygan \cite{Cygan78}) is the unique left-invariant distance on
$\Heis_3$ such that the Cygan distance from $(w_0,w)\in \Heis_3$ to
$(0,0)$ is the {\it Cygan norm} $|(w_0,w)|_{\rm Cyg}= \sqrt{2|w_0|}
$. Note that with the aforementioned change of variables
$(x=\Re\,w,y=\Im\;w, t=2\, \Im\; w_0)$, we do recover the standard
formulation of the Cygan norm on $\RR^3$ (which is, by the way,
equivalent to the Guivarc'h norm introduced much earlier in
\cite{Guivarch70}):
$$
|(x,y,t)|_{\rm Cyg}=\sqrt[4]{(x^2+y^2)^2+t^2}\;.
$$ 
We will denote by $B_{\rm Cyg}(x,r)$ the ball for the Cygan distance
of center $x$ and radius $r$. 

The {\it modified Cygan distance} $d''_{Cyg}$ on $\Heis_3$ is a minor
variation of the Cygan distance, introduced in \cite[\S
4.4]{ParPau11MZ}: it is the unique left-invariant map $d''_{Cyg}$ on
$\Heis_3\times\Heis_3$ such that
$$
d''_{Cyg}((w_0,w),(0,0))=\frac{2|w_0|}{\sqrt{|w^2|+2|w_0|}}\;.
$$ 
Note that $\frac{1}{\sqrt{2}}\;d_{Cyg} \leq d''_{Cyg}\leq d_{Cyg}$.
Though $d''_{Cyg}$ might not be a distance, it is hence close to the
Cygan distance, and it will allow error terms estimates in the
following results.

We consider again the imaginary quadratic number field $K=\QQ(i)$.
The Heisenberg group $\Heis_3$ is the set of real points of an
algebraic group defined over $\QQ$, with group of rational points
equal to $\Heis_3(\QQ)=\Heis_3\cap(K\times K)$. It has also a natural
$\ZZ$-structure, with the  group of integral points equal to
$\Heis_3(\ZZ)=\Heis_3\cap(\OOO_K\times \OOO_K)$.

The following beautiful result (see \cite[Theo.~1.1]{GarNevTay14} for
a more general version), whose tools are in harmonic analysis, solves
the analog in the Heisenberg group of the Gauss circle 
problem. It would be interesting to know if it is still valid for
general imaginary quadratic number field $K$. Note that
\begin{equation}\label{eq:haarvolheisenberg}
\frac{1}{\haarheis(\Heis_3(\ZZ)\bs\Heis_3)}=\frac{2}{|D_K|}
\end{equation} 
(see Equation (24) in \cite{ParPauHeis}, the volume $\Vol_{\Heis_3}$
used in this reference being $\Vol_{\Heis_3}=8\haarheis$), and that
$\haarheis(B_{\rm Cyg} (0,R))=\haarheis(B_{\rm Cyg} (0,1))\,R^4$.

\btheo[Garg-Nevo-Taylor]\label{theo:GaussHeisenberg} As $R\to+\infty$,
we have
$$
\card \big(\Heis_3(\ZZ)\cap B_{\rm Cyg}(0,R)\big)\;=\;
\frac{2\,\haarheis(B_{\rm Cyg} (0,1))}{|D_K|}\;R^4+\bigO(R^2).
$$
\etheo

Any element in $\Heis_3(\QQ)$ may and will be written $(\frac ac,\frac
bc)$, where $a,b,c\in\OOO_K$ are relatively prime (that is, the ideal
of $\OOO_K$ generated by $a,b,c$ is equal to $\OOO_K$), and satisfy
$c\neq 0$ and $2\,\Re(a\overline{c})=|b|^2$. This writing is unique up
to the multiplication of $a,b,c$ by the same element of
$\OOO_K^\times$. In particular, the map $H:\Heis_3(\QQ)\ra\RR$
$$
H\big(\frac ac,\frac bc\big)=|c|
$$ 
is well defined, and is clearly an appropriate height function on 
$\Heis_3(\QQ)$ modulo the action of $\Heis_3(\ZZ)$.

By studying the Diophantine approximation in the Heisenberg group by
its rational points, we understand, as in example (3) at the end of
Subsection \ref{subsec:approxframwork}, taking $X=\Heis_3(\ZZ)\bs
\Heis_3$, $Z=\Heis_3(\ZZ)\bs\Heis_3(\QQ)$, $Y=\Heis_3(\ZZ)\bs
(\Heis_3-\Heis_3(\QQ))$ and the above height function.

The following result is an equidistribution theorem of the set of
rational points in $\Heis_{3}$, analogous to the equidistribution of
Farey fractions in the real field. We denote by $\zeta$ Riemann's zeta
function. Note that the exponent $4$ that appears below is the same as
in the above theorem of Garg-Nevo-Taylor, it is the Hausdorff
dimension of the Cygan distance.

\btheo[Parkkonen-Paulin]\label{theo:equidisHeis} 
As $s\ra+\infty$, we have
$$
\frac{\pi\,|\OOO_K^\times|\,|D_K|^{\frac{3}{2}}\,\zeta_K(3)}
{\zeta(3)}\;s^{-4}\sum_{(\frac ac,\frac
bc)\in\Heis_3(\QQ)\;:\;H\big(\frac ac,\frac bc\big)\leq s}\;
\Delta_{(\frac{a}{c},\frac{\alpha}{c})}\;\weakstar\;\haarheis\,.
$$ 
\etheo

We refer to \cite[Theo.~13]{ParPauHeis} for a version of this
theorem valid for any imaginary quadratic number field, with
additional congruence assumptions, and with error term. The next
result, analogous to the Mertens theorem in the real field, follows from
the version with error term of Theorem \ref{theo:equidisHeis}, by
Equation \eqref{eq:haarvolheisenberg}.

\bcoro[Parkkonen-Paulin]\label{coro:countHeis}  
There exists $\kappa>0$ such that, as $s\ra+\infty$,
$$
\card\;
\Big\{(\frac ac,\frac bc)\in\,\Heis_3(\ZZ)\bs\Heis_3(\QQ)\;:\;
H\big(\frac ac,\frac bc\big)\leq s\Big\}
= \frac{\zeta(3)}
{2\,\pi\,|\OOO_K^\times|\,|D_K|^{\frac{1}{2}}\,\zeta_K(3)}
\,s^4+\bigO(s^{4-\kappa})\,.
$$
\ecoro

We now turn to equidistribution and counting results in the Heisenberg
group of arithmetically defined topological circles, relating them to
Diophantine approximation problems.

Let us consider the Hermitian form of signature $(1,2)$ on $\CC^3$
defined by $$h:(z_0,z_1,z_2)\mapsto -z_0\overline{z_2}-
z_2\overline{z_0} +|z_1|^2\,.
$$ Using homogeneous coordinates in the
complex projective plane $\PP_2(\CC)$, {\it Poincaré's
  hypersphere}\footnote{Actually, Poincaré in \cite{Poincare07} was
  using another Hermitian form with signature $(1,2)$.}  is the
projective isotropic locus of $h$
$$
\H\!\Scal=\{[z_0:z_1:z_2]\in \PP_2(\CC)\;:\; h(z_0,z_1,z_2)=0\}\;,
$$
which is a real-analytic submanifold of $\PP_2(\CC)$ diffeomorphic to
the $3$-sphere $\SSS^3$.  The projective action on $\PP_2(\CC)$ of the
projective special unitary group $\PSU_h$ of $h$ preserves
$\H\!\Scal$. The Alexandrov compactification $\Heis_3\cup \{\infty\}$
of the Heisenberg group $\Heis_3$ identifies with Poincaré's
hypersphere by mapping $(w_0,w)$ to $[w_0:w:1]$ and $\infty$ to
$[1:0:0]$. We identify $\Heis_3$ with its image in $\PP_2(\CC)$,
called {\it Segre's hyperconic}, that we will think of as the
projective model of the Heisenberg group.

As defined by von Staudt,\footnote{Though many references, including
  \cite{Goldman99}, attribute the notion of chains to E.~Cartan, he
  himself attributes them to von Staudt in \cite[footnote
  3)]{Cartan32}.} a {\it chain} in Poincar\' e's hypersphere
$\H\!\Scal$ is an intersection, nonempty and not reduced to a point,
with $\H\!\Scal$ of a complex projective line in $\PP_2(\CC)$. It is
called {\it finite} if it does not contain $\infty=[1:0:0]$. A chain
$C$ separates the complex projective line containing it into two real
discs $D_\pm(C)$, which we endow with their unique Poincar\'e metric
(of constant curvature $-1$) invariant under the stabiliser of $C$ in
$\PSU_h$.

If $\pi:\Heis_3\ra \CC$ is the canonical Lie group morphism
$(w_0,w)\mapsto w$, then the chains are the images, under the elements
of $\PSU_h$, of the {\it vertical chains}, that are the union with
$\{\infty\}$ of the fibers of $\pi$. In particular, the finite chains
are ellipses in (the aforementioned coordinates $(x,y,t)$ of)
$\Heis_3$ whose images under $\pi$ are Euclidean circles in $\CC$. We
refer for instance to \cite[\S 4.3]{Goldman99} for these informations
and more on the chains.

A chain $C$ will be called {\it arithmetic} (over $K=\QQ[i]$) if its
stabiliser $\PSU_h(\OOO_K)_C$ in the arithmetic group $\PSU_h(\OOO_K)=
\PSU_h\cap \operatorname{PGL}_3(\OOO_K)$ has a dense orbit in $C$. It
then turns out that $\PSU_h(\OOO_K)_C$ acts discretely with finite
covolume on $D_\pm(C)$, and we denote by $\covol(C)$ the (common)
volume of $\PSU_h(\OOO_K)_C\bs D_\pm(C)$. Note that the stabiliser
$\PSU_h(\OOO_K)_\infty$ of $[1:0:0]$ in $\PSU_h(\OOO_K)$ preserves the
$d''_{\rm Cyg}$-diameters of the chains. The picture below 
shows part of an orbit of arithmetic chains under the arithmetic lattice
$\PSU_h(\OOO_K)$.

\begin{center}
\includegraphics[width=9cm]{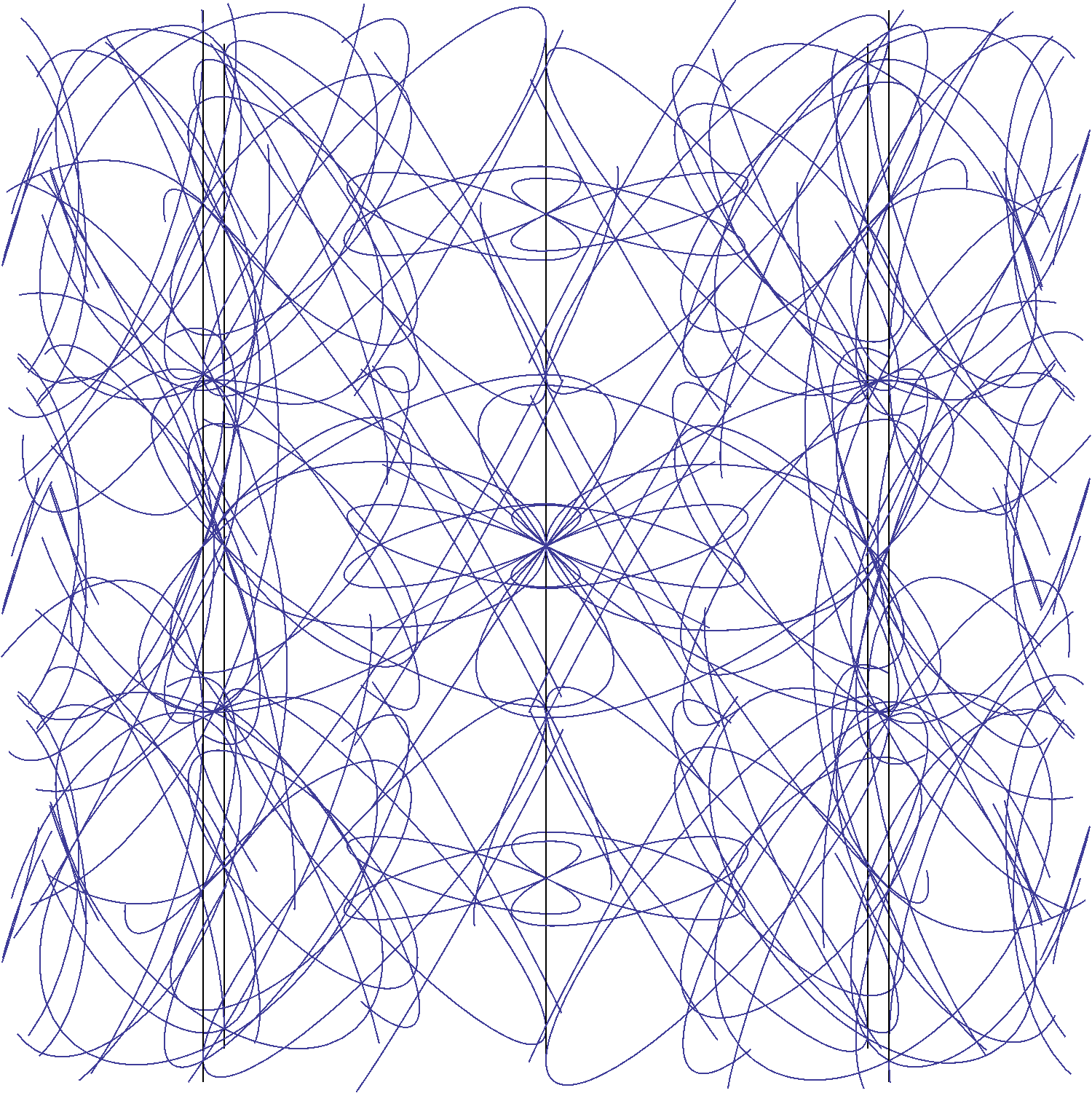}
\end{center}

\btheo[Parkkonen-Paulin]\label{theo:countchain} Let $C_0$ be an
arithmetic chain over $K$ in the hypersphere $\H\!\Scal$. Then there
exists a constant $\kappa>0$ such that, as $\epsilon>0$ tends to $0$,
the number of chains modulo $\PSU_h(\OOO_K)_\infty$ in the
$\PSU_h(\OOO_K)$-orbit of $C_0$ with $d''_{\rm Cyg}$-diameter at least
$\epsilon$ is equal to
$$
\frac{512\,\zeta(3)\,\covol(C_0)}
{|\OOO_K^\times|\,|D_K|^{\frac{3}{2}}\,\zeta_K(3)\,n_{0}}
\;\epsilon^{-4}\big(1+\bigO(\epsilon^\kappa)\big)\;,
$$
where $n_{0}$ is the order of the pointwise stabiliser of $C_0$ in
$\PSU_h(\OOO_K)$.  
\etheo

We refer to \cite[Theo.~19]{ParPauHeis} for a version of this theorem
valid for any imaginary quadratic number field (the pictures below
represent two views of a part of an orbit of an arithmetic chain when
$K=\QQ(i\sqrt{2})$)
and with additional congruence assumptions.

\begin{center}
\includegraphics[width=7cm]{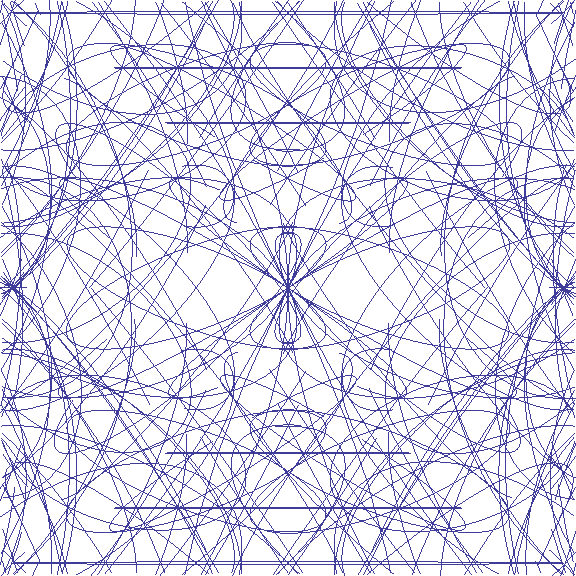}\hfill
\includegraphics[width=7cm]{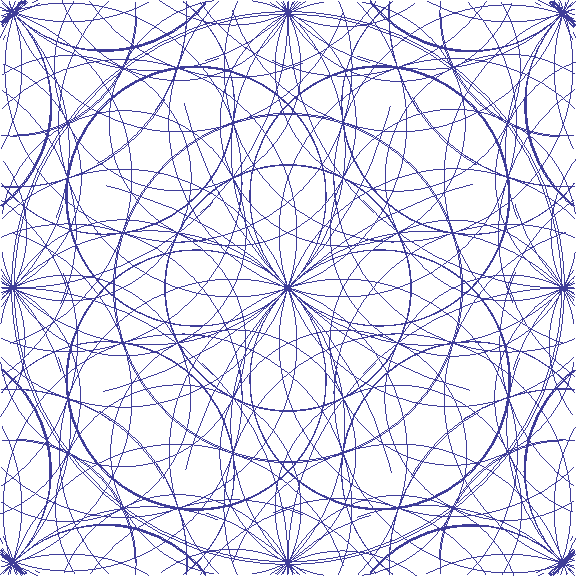}
\end{center}

Given a complex projective line $L$ in $\PP_2(\CC)$, there is a unique
order $2$ complex projective map with fixed point set $L$, called the
{\it reflexion} in $L$. Given a finite chain $C$, contained in the
complex projective line $L(C)$, the {\it center} of $C$ (see for
instance \cite[4.3.3]{Goldman99}), denoted by $\operatorname{cen}(C)
\in \H\!\Scal-\{\infty\}=\Heis_3$, is the image of $\infty=[1:0:0]$
under the reflexion in $L(C)$. We also prove in
\cite[Theo.~20]{ParPauHeis} (a version valid for any imaginary
quadratic number field and allowing additional congruence assumptions
of) the following result saying that the centers of the finite
arithmetic chains in a given $\PSU_h(\OOO_K)$-orbit of a given
arithmetic chain $C_0$ with $d''_{\rm Cyg}$-diameter at least
$\epsilon$ equidistribute in the Heisenberg group towards the
normalised Haar measure $\frac{2}{|D_K|} \,\haarheis$. 

\btheo[Parkkonen-Paulin]\label{theo:equidischain} Let $C_0$ be an
arithmetic chain over $K$. As $\epsilon>0$ tends to $0$, we have
$$
\frac{n_{0}\,(1+2\,\delta_{D_K,-3})
  \,|D_K|^{\frac{5}{2}}\,\zeta_K(3)}
{128\,\zeta(3)\,\covol(C_0)}\;\epsilon^{4}
\sum_{C\in \PSU_h(\OOO_K)\cdot C_0\;:\;\diam_{d''_{\rm Cyg}}\;C\geq \epsilon}\;
\Delta_{\operatorname{cen}(C)}\;\weakstar\;\haarheis\,.
$$ 
\etheo
The above result fits into the framework of Subsection
\ref{subsec:approxframwork} with
$$
X=Y=\Heis_3(\ZZ)\bs \Heis_3,\;\;\; Z=\{\ga \operatorname{cen}(C_0)
\;:\;\ga\in\PSU_h(\OOO_K), \;\infty\notin \ga C_0\}
$$ 
and $H: \operatorname{cen}(\ga C_0) \mapsto
\diam_{d''_{\rm Cyg}}(\ga C_0)$.

\medskip This equidistribution phenomenon can be understood by looking
at the following picture, representing a different view of the same
orbit of arithmetic chains as the one above Theorem
\ref{theo:countchain}.

\begin{center}
\includegraphics[width=11cm]{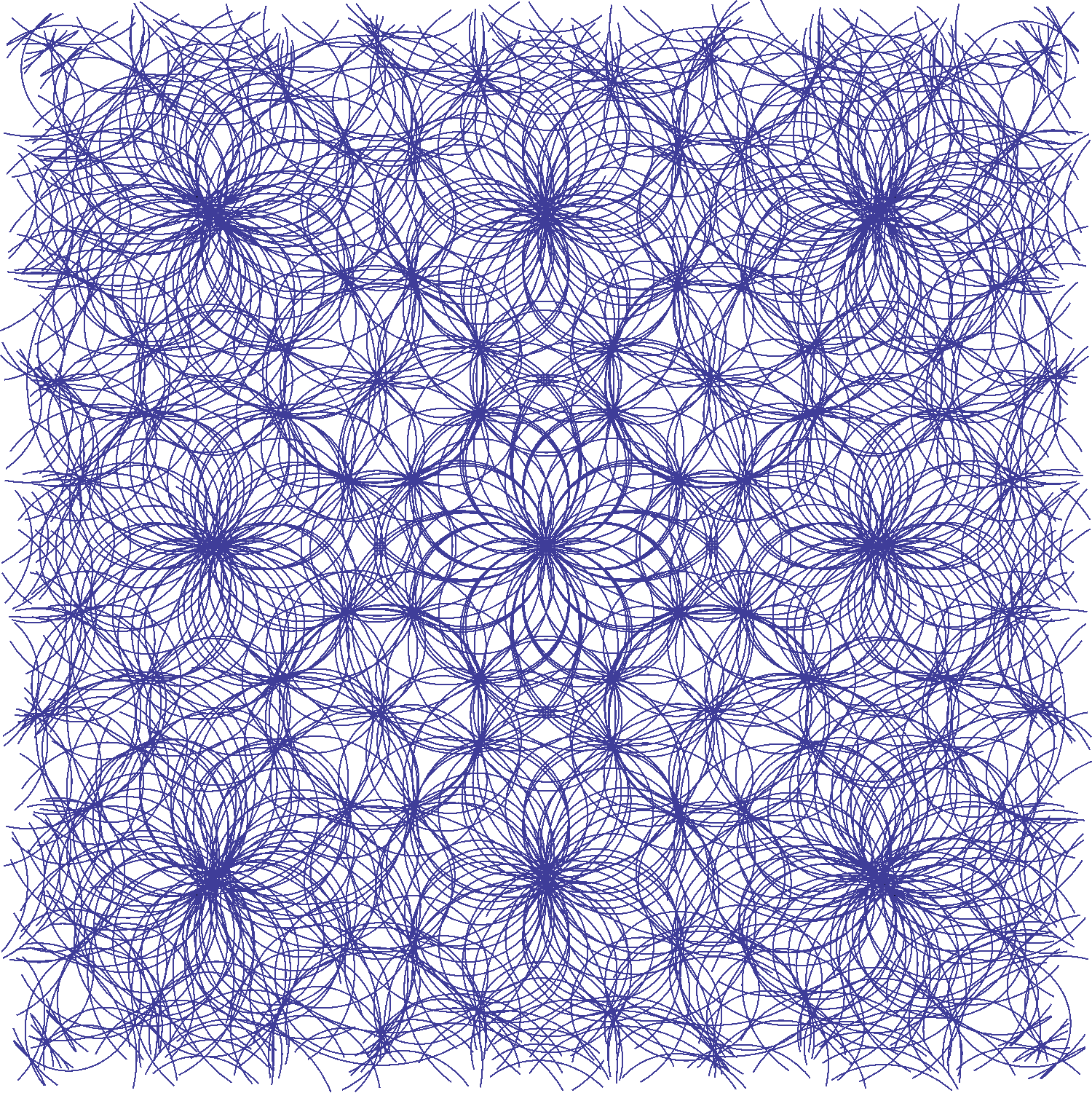}
\end{center}

\section{Measures in negative curvature}
\label{sect:measures}

\subsection{A classical link between basic Diophantine 
approximation and hyperbolic geometry}
\label{subsec:classiclink}

Let us briefly explain a well-known link between hyperbolic geometry
and Diophantine approximation problems, which goes back to Gauss and
Ford (see also \cite{Series85a}). This will start to explain why the
proofs of Theorems \ref{theo:Lagrangequadirrat},
\ref{theo:equidistrace}, \ref{theo:relatcomplcount},
\ref{theo:equidistribR}, \ref{theo:equidistribC},
\ref{theo:equidisHeis}, \ref{theo:countchain} and
\ref{theo:equidischain} all use real or complex hyperbolic geometry.

For $n\geq 2$, let $\HH^n_\RR$ be the upper halfspace model of the
real hyperbolic $n$-space, with underlying manifold
$\{x=(x_1,\dots,x_n) \in\RR^n\;:\; x_n>0\}$ and Riemannian metric
$\frac{|dx|^2}{x_n^2}$.

Within the framework of example (1) of Subsection
\ref{subsec:approxframwork}, let $\partial_\infty\HH^2_\RR= \RR\cup
\{\infty\}$ be the boundary at infinity of $\HH^2_\RR$, $Y=\RR$ and
$Z= \QQ= \Ga\cdot\infty- \{\infty\}$, where $\Ga=\PSL_2(\ZZ)$ is the
well-known {\it modular group}, which is a nonuniform lattice in
$\PSL_2(\RR)= \Isom_0(\HH^2_\RR)$, acting by homographies on $\RR$.

\medskip\noindent
\begin{minipage}{8.4cm}
  ~~~ Let $\psi:\QQ\cup\{\infty\}\ra [1,+\infty[$ be a map, let
  $\H_\infty=\{z\in\CC\;:\Im\;z\geq \psi(\infty)\}$ and, for every
  $\frac pq\in\QQ$, let $\H_{\frac pq}$ be the closed Euclidean disc
  of center $\frac pq +\frac{1}{2\, \psi(\frac pq)q^2}i$ and radius
  $\frac{1}{2\,\psi(\frac pq)q^2}$, with its tangency point to the
  horizontal axis removed. Then $(\H_x)_{x\in\QQ\cup\{\infty\}}$ is a
  family of horodiscs in $\HH^2_\RR$, with pairwise disjoint
  interiors, centred at the parabolic fixed points of $\Ga$. If $\psi$
  is constant, then this family is $\Ga$-equivariant, and if $\psi=1$,
  then it is maximal.
\end{minipage}
\begin{minipage}{6.5cm}
\begin{center}
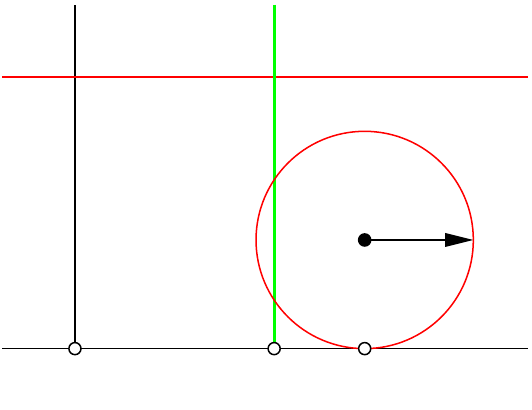
\end{center}
\end{minipage}

\medskip
A link between Diophantine approximation of real numbers by rational
ones and hyperbolic geometry is the following one: for every
$\xi\in\RR$, we have
$$
\big|\xi-\frac pq\,\big|\leq \frac{1}{2\,\psi(\frac pq)q^2}
$$
if and only if the geodesic line $L_\xi$ in $\HH^2_\RR$ from
$\infty\in\partial_\infty\HH^2_\RR$ to $\xi\in \partial_\infty
\HH^2_\RR$ meets the horodisc $\H_{\frac pq}$.  Many Diophantine
approximation properties of $\xi$ may be explained by the behaviour of
the image of the geodesic $L_\xi$ in the modular curve $\Ga\bs
\HH^2_\RR$. For instance, the coefficients of the continued fraction
expansion of $\xi\in\RR-\QQ$ are bounded if and only if the positive
subray of $L_\xi$ has a bounded image in $\Ga\bs\HH^2_\RR$ (see for
instance \cite{Series85a}).

\medskip
Hence, it is useful for arithmetic applications to study the dynamical
and ergodic properties of the geodesic flows in negative curvature,
and we develop these topics in the following two sections.

\subsection{Negative curvature background}
\label{subsec:negacurv}

We refer for instance to \cite{BriHae99,Roblin03,PauPolSha} for
definitions, proofs and complements concerning this subsection.

Let $\wt M$ (for instance $\HH^2_\RR$) be a complete simply connected
(smooth) Riemannian manifold with (dimension at least $2$ and) pinched
negative sectional curvature $-b^2\le K\le -1$, and let $x_0\in\wt M$
be a fixed basepoint.  Let $\Ga$ (for instance $\PSL_2(\ZZ)$) be a
nonelementary (not virtually nilpotent) discrete group of isometries
of $\wt M$, and let $M$ be the quotient Riemannian orbifold $\Ga\bs\wt
M$. See \cite{Roblin03, BroParPau15} to relax the pinching and
manifold assumptions on $\wt M$.

We denote by $\partial_{\infty}\wt M$ the boundary at infinity of $\wt
M$, that is the quotient space of the space of geodesic rays
$\rho:[0,+\infty[\;\ra\wt M$ in $\wt M$, two of them being equivalent
if the Hausdorff distance between their images is finite. The class of a
geodesic ray is called its {\it point at infinity}.

The {\it Busemann cocycle} of $\wt M$ is the map $\beta: \wt
M\times\wt M\times\partial_{\infty} \wt M\to\RR$ defined by $$
(x,y,\xi)\mapsto \beta_{\xi}(x,y)=
\lim_{t\to+\infty}d(x,\rho(t))-d(y,\rho(t))\;,
$$
where $\rho$ is any geodesic ray with point at infinity $\xi$. The
above limit exists and is independent of $\rho$.  The {\it horosphere}
with center $\xi\in\partial_{\infty}\wt M$ through $x\in\wt M$ is
$\{y\in\wt M\;:\; \beta_{\xi}(x,y)=0\}$, and $\{y\in\wt M\;:\;
\beta_{\xi}(x,y)\leq 0\}$ is the {\it horoball} centered at $\xi$
bounded by this horosphere. For instance, in $\HH^n_\RR$, the
horoballs are either the subspaces $\{(x_1,\dots,x_n)\in \HH^n_\RR
\;:\; x_n\geq a\}$ for $a>0$ or the closed Euclidean balls contained
in the closure of $\HH^n_\RR$, tangent to the horizontal hyperplane,
minus the tangency point.

Let $\wt M\cup \partial_{\infty}\wt M$ be the geometric
compactification of $\wt M$, and let $\Lambda\Ga=\overline{\Ga
  x_0}-\Ga x_0$ be the {\it limit set} of $\Ga$. Let $\pi:T^1\wt M\ra
\wt M$ be the unit tangent bundle of $\wt M$.

\medskip\noindent
\begin{minipage}{9.4cm}
  ~~~ For every $v\in T^1\wt M$, let $v_-\in\partial_\infty\wt M$ and
$v_+\in\partial_\infty\wt M$, respectively, be the endpoints at
$-\infty$ and $+\infty$ of the geodesic line defined by $v$.  Let
$\partial_\infty^2\wt M$ be the subset of $\partial_\infty\wt
M\times\partial_\infty\wt M$ which consists of pairs of distinct
points at infinity of $\wt M$.  {\em Hopf's parametrisation} of
$T^1\wt M$ is the homeomorphism which identifies $T^1\wt M$ with
$\partial_\infty^2\wt M\times\RR$, by the map $v\mapsto(v_-,v_+,s)$,
where $s$ is the signed distance of the closest point to $x_0$ on the
geodesic line defined by $v$ to $\pi(v)$.  
\end{minipage}
\begin{minipage}{5.5cm}
\begin{center}
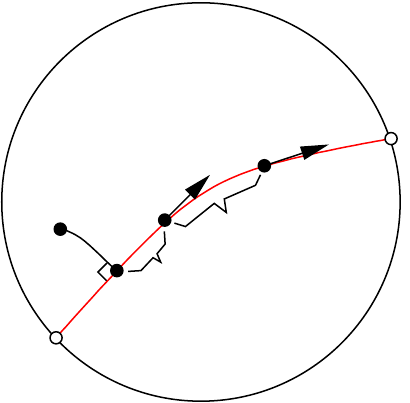
\end{center}
\end{minipage}

\medskip The geodesic flow is the smooth flow $(\flow{t})_{t\in\RR}$
on $T^1\wt M$ defined, in Hopf's coordinates, by $\flow{t}: (v_-, v_+,s)
\mapsto (v_-,v_+,s+t)$.

\medskip 
The {\em strong stable manifold} of $v\in T^1\wt M$ is 
$$
W^{+}(v)=\{v'\in T^1\wt M\;:\;
\lim_{t\to+\infty}d(\pi(\flow{t}v),\pi(\flow{t}v'))=0\},  
$$
and  its {\em strong unstable manifold} is 
$$
W^{-}(v)=\{v'\in T^1\wt M\;:\;
\lim_{t\to-\infty}d(\pi(\flow{t}v),\pi(\flow{t}v'))= 0\}. 
$$
The {\em stable/unstable manifold} of $v\in T^1\wt M$ is
$W^{0\pm}(v)=\bigcup_{t\in\RR}\flow t W^{\pm}(v)$, which consists of
the elements $v'\in T^1\wt M$ with $v'_\pm=v_\pm$. The map $(t,w)$
from $\RR \times W^{\pm}(v)$ to $W^{0\pm}(v)$ is a homeomorphism. The
projections $\pi(W^{+}(v))$ and $\pi(W^{-}(v))$ in $\wt M$ of the
strong stable and strong unstable manifolds of $v$ are the horospheres
through $\pi(v)$ centered at $v_{+}$ and $v_{-}$, respectively (see
the picture below on the left hand side). These horospheres bound
horoballs denoted by $\operatorname{HB}^\pm(v)$. The strong stable
manifolds and strong unstable manifolds are the (smooth) leaves of
topological foliations in $T^1\wt M$ that are invariant under the
geodesic flow and the group of isometries of $\wt M$, denoted by
$\W^{+}$ and $\W^{-}$ respectively.

\begin{center}
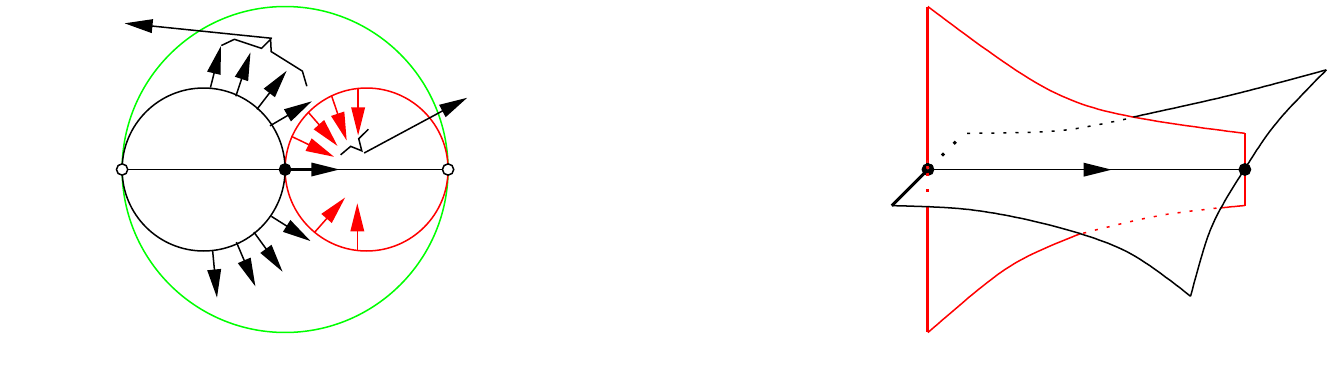
\end{center}

For every $v\in T^1\wt M$, let $d_{W^{-}(v)}$ and
$d_{W^{+}(v)}$ be {\it Hamenst\"adt's distances} on the strong
unstable and strong stable leaf of $v$, defined as follows (see
\cite[Appendix]{HerPau97} and compare with \cite{Hamenstadt89}): for all
$w,z\in W^{\pm}(v)$, let
$$
d_{W^{\pm}(v)}(w,z) = \lim_{t\ra+\infty} 
e^{\frac{1}{2}d(\pi(\phi_{\mp t}w),\;\pi(\phi_{\mp t}z))-t}\;.
$$
The above limit exists, and Hamenst\"adt's distance induces the
original topology on $W^{\pm}(v)$, though it has fractal properties:
the Hausdorff dimension of $d_{W^{+}(v)}$ is in general bigger than
the topological dimension of $W^{\pm}(v)$. For all $w,z\in W^{\pm}(v)$
and $t\in\RR$, and for every isometry $\ga$ of $\wt M$, we have
$$
d_{W^{\pm}(\ga v)}(\ga w,\ga z)= d_{W^{\pm}(v)}(w,z)\quad{\rm and}\quad
d_{W^{\pm} (\flow tv)}(\flow tw,\flow tz)=e^{\mp t}d_{W^{\pm}(v)}(w,z)\;.
$$
These dilation/contraction properties of Hamenst\"adt's distances
under the geodesic flow are a strengthening of the Anosov property of
the geodesic flow (see the above picture on the right hand side).

\medskip
The {\it Poincar\'e series} of $\Ga$ is the map $P_{\Ga}:
\RR\ra\mathopen{]}0,+\infty\mathclose{]}$ defined by
$$
P_{\Ga}(s)=\sum_{\ga\in\Ga} \;\; e^{-s\,d(\ga x_0,\,x_0)}\;.
$$
The {\it critical exponent} of $\Ga$ is
$$
\delta_{\Ga}=\lim_{n\ra +\infty}\;\frac{1}{n}\ln
\;\card\{\ga\in\Ga,\; d(x_0,\,\ga x_0)\leq n\}\;.
$$
The above limit exists, is independent of $x_0$, we have
$\delta_{\Ga}\in \mathopen{]}0,+\infty\mathclose{[}$ and the
Poincar\'e series $P_\Ga(s)$ of $\Ga$ converges if $s>\delta_\Ga$ and
diverges if $s<\delta_\Ga$, see for instance \cite{Roblin02, Roblin03},
as well as \cite{PauPolSha} for versions with potential).

\subsection{The various measures}
\label{subsec:variousmeasures}

We refer for instance to \cite{Roblin03,ParPau14ETDS,PauPolSha,
  BroParPau15} for definitions, proofs and complements concerning this
subsection. We introduce here the various measures that will be useful
for our ergodic study in Section \ref{sect:geomequidistribcounting}

\medskip A {\it Patterson-Sullivan measure} is a family $(\mu_x)_{x\in
  \wt M}$ of finite nonzero measures on $\partial_\infty \wt M$, whose
support is the limit set $\Lambda\Ga$ of $\Ga$, such that, for all
$\ga\in\Ga$, $x,y\in \wt M$ and $\xi\in\partial_\infty\wt M$,
$$
\ga_*\mu_x=\mu_{\ga x}\quad{\rm and}\quad
d\mu_x(\xi)=e^{-\delta_\Ga\,\beta_\xi(x,\,y)}\;d\mu_y(\xi)\,.
$$
Such a family exists, and if $P_\Ga(\delta_\Ga)=+\infty$, then, for all
$x\in \wt M$,
$$
\mu_x=\lim_{s\ra {\delta_\Ga}^+}\frac{1}{P_\Ga(s)}\;
\sum_{\ga\in\Ga} \;\; e^{-s\,d(x,\ga x_0)}\;\Delta_{\ga x_0}
$$
for the weak-star convergence of measures (see for instance \cite{Roblin03}).

\medskip Let $C$ be a nonempty proper closed convex subset of $\wt M$,
with stabiliser $\Ga_C$ in $\Ga$, such that the family $(\ga
C)_{\ga\in \Ga/\Ga_C}$ of subsets of $\wt M$ is locally finite.

\smallskip\noindent
\begin{minipage}{10.4cm}
  ~~~ The {\em inner} (respectively {\em outer}) {\em unit normal
    bundle} $\normalin C$ (respectively $\normalout C$) of $C$ is the
  topological submanifold of $T^1\wt M$ consisting of the unit tangent
  vectors $v\in T^1\wt M$ such that $\pi(v)\in\partial C$, $v$ is
  orthogonal to a contact hyperplane to $C$ and points towards
  (respectively away from) $C$ (see \cite{ParPau14ETDS} for more
  precisions, and note that the boundary $\partial C$ of $C$ is not
  necessarily $C^1$, hence may have more than one contact hyperplane
  at some point, and that it is not necessarily true that $\exp(tv)$
  belongs (respectively does not belong) to $C$ for $t>0$ small
  enough).
\end{minipage}
\begin{minipage}{4.5cm}
\begin{center}
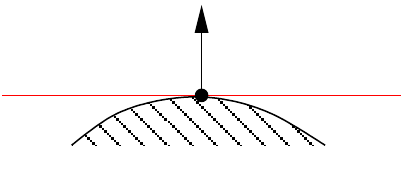
\end{center}
\end{minipage}

\smallskip The endpoint map $w\mapsto w_\pm$ from $\normalpm C$ into
$\partial_\infty \wt M- \partial_\infty C$ is a homeomorphism. Using
this homeomorphism, we defined in \cite{ParPau14ETDS} the {\it outer}
(respectively {\it inner}) {\it skinning measure} on $T^1\wt M$ as the
measure $\wt\sigma_{C}^+$ (respectively $\wt\sigma_{C}^-\,$) with
support in $\normalout C$ (respectively $\normalin C$) given by, for
all $w\in\normalpm C$,
\begin{equation}\label{eq:defiskinmeas}
d\wt\sigma_{C}^\pm(w)  =  
e^{-\delta_{\Ga}\beta_{w_{\pm}}(\pi(w),\,x_{0})}\,d\mu_{x_{0}}(w_{\pm})\,.
\end{equation}
This measure is independent of $x_0$ and satisfies $\ga_*\wt
\sigma_{C}^\pm=\wt\sigma_{\ga C}^\pm$ for every $\ga\in\Ga$.  Hence
the measure $\sum_{\ga\in\Ga/\Ga_C} \wt\sigma_{\ga C}^\pm(w)$ is a
$\Ga$-invariant locally finite measure on $T^1\wt M$, therefore
defining (see for instance \cite[\S 2.4]{ParPau13b} for details) a
(locally finite) measure on $T^1M=\Ga\bs T^1\wt M$, called the {\em
  outer} (respectively {\it inner}) {\it skinning measure} of $C$ on
$T^1M$, and denoted by $\sigma_{C}^+$ (respectively $\sigma_{C}^-\,$).

Note that the measure $\wt\sigma_{\operatorname{HB}^-(w)}^+$
(respectively $\wt\sigma_{\operatorname{HB}^+(w)}^-\,$) coincides with
the Margulis measure (see for instance \cite{Margulis04a, Roblin03})
on the strong unstable leaf $W^-(w)$ (respectively strong stable leaf
$W^+(w)\,$), for every $w\in T\wt M$.

When $\wt M$ has constant curvature and $\Ga$ is geometrically finite,
when $C$ is an ball, horoball or totally geodesic submanifold, the
(outer) skinning measure of $C$ has been introduced by Oh and Shah
\cite{OhSha12,OhSha13}, who coined the term, with beautiful
applications to circle packings, see also \cite[Lemma 4.3]{HerPau10}
for a closely related measure. The terminology comes from McMullen's
proof of the contraction of the skinning map (capturing boundary
information for surface subgroups of $3$-manifold groups) introduced
by Thurston to prove his hyperbolisation theorem.

\medskip The {\em Bowen-Margulis measure} on $T^1\wt M$ (associated to
a given Patterson-Sullivan measure) is the measure $\wt m_{\rm BM}$ on
$T^1\wt M$ given by the density
\begin{equation}\label{eq:defigibbs}
d\wt m_{\rm BM}(v)=
e^{-\delta_\Ga(\beta_{v_-}(\pi(v),\,x_0)\,+\,\beta_{v_+}(\pi(v),\,x_0))}\;
d\mu_{x_0}(v_-)\,d\mu_{x_0}(v_+)\,ds
\end{equation} 
in Hopf's parametrisation. The Bowen-Margulis measure $\wt m_{\rm BM}$
is independent of $x_0$, and it is invariant under the actions of the
group $\Ga$ and of the geodesic flow. Thus (see for instance \cite[\S
2.4]{ParPau13b}), it defines a measure $m_{\rm BM}$ on $T^1M$ which is
invariant under the quotient geodesic flow, called the {\em
  Bowen-Margulis measure} on $T^1M$. We refer for instance to
\cite{Ledrappier95b,PauPolSha} for the extensions to the case with
potential of the Patterson-Sullivan and Bowen-Margulis measures (the
latter becomes the Gibbs measure), and to \cite{ParPau13b} for the
extensions to the case with potential of the skinning measures.

\medskip 
Let $C$ be a nonempty proper closed convex subset of $\wt M$.

\smallskip\noindent
\begin{minipage}{7.9cm}
 Let 
$$U_C^\pm=\{v\in T^1\wt M\;:\;v_\pm\notin \partial_\infty
  C\}\;,
$$
  and let $f_C: U_C^\pm\ra \normalpm C$ be the map sending $v$
  to the unique $w\in \normalpm C$ such that $w_\pm=v_\pm$. It is a
  topological fibration, whose fiber over $w\in \normalpm C$ is its
  stable/unstable leaf $W^{0\pm}(w)$.
\end{minipage}
\begin{minipage}{7cm}
\begin{center}
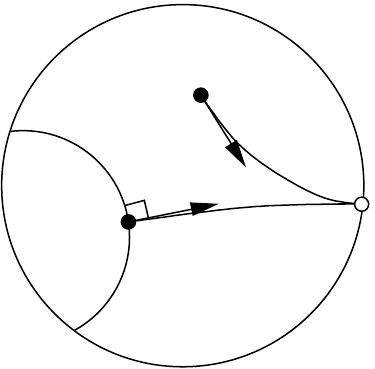
\end{center}
\end{minipage}

\smallskip 
Given $w\in T^1\wt M$, since $\normalmp \operatorname{HB}^{\pm}(v)
=W^{\pm}(v)$ and using the homeomorphism $(t,v)\mapsto \flow t v$ from
$\RR\times W^{\pm}(w)$ to $W^{0\pm}(w)$, we define a measure on
$T^1\wt M$ with support contained in $W^{0\pm}(w)$ by, with $v\in
W^{\pm}(w)$ and $t\in\RR$,
$$
d\mu_{W^{0\pm}(w)}(\flow t v)=dt\;d\wt\sigma^\mp_{\operatorname{HB}^{\pm}(w)}(v) \;.
$$
Note that the measure $d\mu_{W^{0\pm}(w)}$ depends on $W^{\pm}(w)$,
not only on $W^{0\pm}(w)$.

The following result (see \cite[Prop.~8]{ParPau14ETDS}) says that the
Bowen–Margulis measure disintegrates over the skinning measure of
$C$. When $\wt M$ has constant curvature and $\Ga$ is torsionfree and
for special convex sets $C$, this result is implicit in
\cite{OhSha13}.

\bprop[Parkkonen-Paulin]\label{prop:disintegration} The restriction to
$U_C^\pm$ of the Bowen-Margulis measure $\wt m_{\rm BM}$ disintegrates
by the fibration $f_C^\pm:U_C^\pm\ra \normalpm C$, over the skinning
measure $\wt\sigma_{C}^\pm$ of $C$, with conditional measure
$\mu_{W^{0\pm}(w)}$ on the fiber $W^{0\pm}(w)$ of $f_C^{\pm}$ over
$w\in \normalpm C$: for $v\in U_C^\pm$,
$$
d\wt m_{\rm BM}(v)=
\int_{w\in \normalpm C}\;
d\mu_{W^{0\pm}(w)}(v)\;d\wt\sigma_{C}^\pm(w)\;.
$$
\eprop

We summarize in the following statement the ergodic properties of the
Bowen-Margulis measure that we will use in Section
\ref{sect:geomequidistribcounting}.

\btheo \label{theo:requisitgeodflow} Assume that $m_{\rm BM}$ is finite.
\begin{enumerate}
\item[(1)] {\bf (Patterson, Sullivan, Roblin)} We have
  $P_\Ga(\delta_\Ga)=+\infty$ and the Patterson-Sullivan measure is
  unique up to a multiplicative constant; hence the Bowen-Margulis
  measure $m_{\rm BM}$ is uniquely defined, up to a multiplicative
  constant.
\item[(2)] {\bf (Bowen, Margulis, Otal-Peigné)} When normalised to be
  a probability measure, the Bowen-Margulis measure on $T^1M$ is the
  unique measure of maximal entropy of the geodesic flow.
\item[(3)]{\bf (Babillot)} If the set of lengths of closed geodecics
  in $M$ generates a dense subgroup of $\RR$, then $m_{\rm BM}$ is
  mixing under the geodesic flow.
\item[(4)]{\bf (Kleinbock-Margulis, Clozel)} If $\wt M$ is a symmetric
  space and $\Ga$ an arithmetic lattice, then there exist $c,\kappa>0$
  and $\ell\in\NN$ such that for all $\phi,\psi\in \C_c^\ell(T^1M)$
  and all $t\in\RR$, we have
\end{enumerate}
$$
\Big|\int_{T^1M} \phi\circ g^{-t}\;\psi\;dm_{\rm BM}-
\frac{1}{\|m_{\rm BM}\|}\;\int_{T^1M}
\phi\;dm_{\rm BM}\int_{T^1M} \psi\;dm_{\rm BM}\;\Big|\leq
c\,e^{-\kappa |t|}\;\|\psi\|_\ell\;\|\phi\|_\ell\;.
$$
\etheo    

Here are a few comments on these results. The Bowen-Margulis measure
$m_{\rm BM}$ is finite for instance when $M$ is compact, or when $\Ga$
is geometrically finite and its critical exponent is strictly bigger
than the critical exponents of its parabolic subgroups (as it is the
case when $M$ is locally symmetric), by \cite{DalOtaPei00}.

For the second assertion, we refer to \cite{Margulis70} and
\cite{Bowen72a} when $M$ is compact, and to \cite{OtaPei04} under the
weaker assumption that $m_{\rm BM}$ is finite. We refer to
\cite{Roblin03} for a proof of the first assertion, to
\cite[Thm.~1]{Babillot02b} for the third one, and to \cite{PauPolSha}
for the extensions of the first three assertions to the case with
potential. The assumption of Assertion (3), called {\it
  non-arithmeticity of the length spectrum} holds for instance, by
\cite{Dalbo99,Dalbo00}, when $M$ is locally symmetric or
$2$-dimensional, or when $\Ga$ contains a parabolic element, or when
$\Lambda\Ga$ is not totally disconnected.

In Assertion (4), called the {\it exponential decay of correlation}
property of the Bowen-Margulis measure, we denote by $\|\cdot\|_\ell$
the Sobolev norm of regularity $\ell$. We refer to \cite{KleMar96} for
a proof of this last assertion, with the help of \cite[Theorem
3.1]{Clozel03} to check its spectral gap property and of
\cite[Lemma~3.1]{KleMar99} to deal with finite cover problems.  Note
that the spectral gap property has been checked by \cite{MohOh14} if
$\wt M=\HH^n_\RR$ and $\Ga$ is only assumed to be geometrically finite
with $\delta_\Ga\geq n-2$ if $n\geq 3$ and $\delta_\Ga\geq
\frac{1}{2}$ if $n=2$, thus providing the first infinite volume
examples for which Assertion (4) holds.

When $M$ is locally symmetric with finite volume, the Bowen-Margulis
measure $\wt m_{\rm BM}$ coincides, up to a multiplicative constant,
with the Liouville measure, that is the Riemannian measure
$\Vol_{T^1\wt M}$ of Sasaki's metric on $T^1\wt M$.  See for instance
\cite[\S 7]{ParPauRev} when $M$ is real hyperbolic. In particular, the
measure is finite in this case. More precisely (see \cite[\S
7]{ParPauRev} and \cite[\S 6]{ParPau13b}), if $\wt M=\HH^n_\RR$ and if
$M$ has finite volume, normalizing the Patterson-Sullivan measure so
that its total mass is the volume of the $(n-1)$-sphere with its
standard spherical metric (so that $\|\mu_{x_0}\|=\Vol(\SSS^{n-1})$),
we have, denoting by $\Vol_N$ the Riemannian volume of any Riemannian
manifold $N$,
\begin{itemize}
\item $\wt m_{ \rm BM}=2^{n-1}\Vol_{T^1\wt M}$;
\item  if $C$ is a horoball, then $\wt\sigma^\pm_{D}= 2^{n-1}
  \Vol_{\normalpm D}$;
\item if $C$ is totally geodesic, then
  $\wt\sigma^+_D=\wt\sigma^-_D= \Vol_{\normalpm D}$.
\end{itemize}

\section{Geometric equidistribution and counting}
\label{sect:geomequidistribcounting}

In this section, we will link the Diophantine approximation problems
of Subsections \ref{subsec:DiophapproxR}, \ref{subsec:equidistribRC}
and \ref{subsec:equidistribHeisenberg} to  general geometric
equidistribution and counting problems on the common perpendiculars
between two locally convex subsets in a negatively curved Riemannian
orbifold.

Assume for simplicity in the introduction of this Section
\ref{sect:geomequidistribcounting} (thus avoiding problems of
regularity, multiplicities and finiteness), that $N$ is a compact
negatively curved Riemannian manifold, and that $D^-$ and $D^+$ are
proper nonempty disjoint closed locally convex subsets of $M$ with
smooth boundaries. A {\it common perpendicular} from $D^-$ to $D^+$ is
a locally geodesic path in $N$ starting perpendicularly from $D^-$ and
arriving perpendicularly to $D^+$. There is exactly one such common
perpendicular in every homotopy class of paths starting from $D^-$ and
ending in $D^+$, where during the homotopy the origin of the path
remains in $D^-$ and its terminal point remains in $D^+$. In
particular, there are at most countably many such common
perpendiculars, and at most finitely many when their length is
bounded.

Even when $N$ is a closed hyperbolic surface and $D^-, D^+$ are simple
closed geodesic (see the picture below), the result (see Equation
\eqref{eqcaseclosedgeod}) was not known before appearing in
\cite{ParPau13b}.

\begin{center}
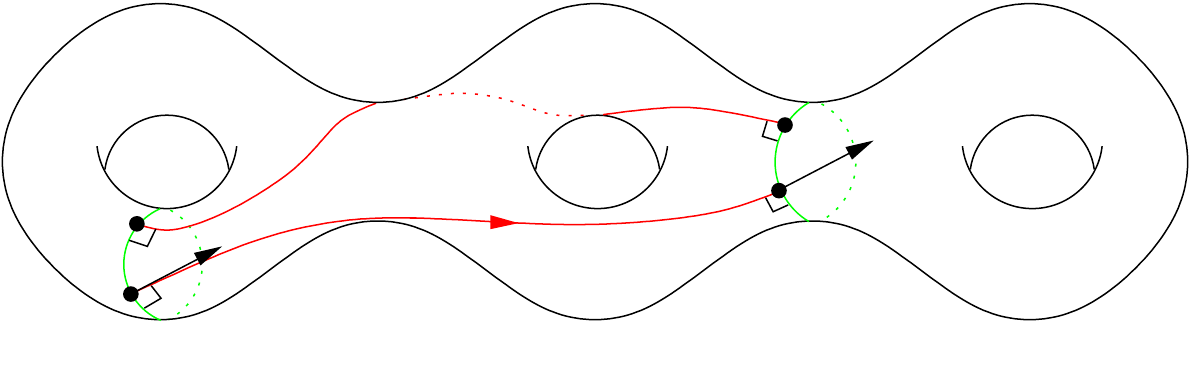
\end{center}

We give in Subsection \ref{subsect:commonperp} (refering to
\cite{ParPau13b} for complete statements and proofs) an asymptotic
formula as $t\ra+\infty$ for the number of common perpendiculars of
length at most $t$ from $D^-$ to $D^+$, and an equidistribution result
as $t\ra+\infty$ of the initial and terminal tangent vectors
$v^-_\alpha$ and $v^+_\alpha$ of these common perpendiculars $\alpha$
in the outer and inner unit normal bundles of $D^-$ and $D^+$,
respectively. Although we do use Margulis's mixing ideas, major new
techniques needed to be developed to treat the problem in the
generality considered in \cite{ParPau13b}, some of them we will
indicate in Subsection \ref{subsect:commonperp}.

\medskip Here is a striking corollary of Theorem \ref{theo:maincount}
in a very different context, that apparently does not involve negative
curvature dynamics or geometry. Let $\Ga$ be a geometrically finite
discrete subgroup of $\PSL_{2}(\CC)$ (acting by homographies on
$\PP_1(\CC)=\CC\cup\{\infty\}$). Assume that $\Ga$ does not contain a
quasifuchsian subgroup with index at most $2$, and that its limit set
$\Lambda\Ga$ is bounded and not totally disconnected in $\CC$. These
assumptions are only here to ensure that the domain of discontinuity
$\Omega\Ga= (\CC\cup\{\infty\})-\Lambda\Ga$ of $\Ga$ has infinitely
many connected components (only one of them unbounded). The following
result gives a precise asymptotic as $\epsilon$ tends to $0$ on the
counting function of the (finite) number of these connected components
whose diameter are at least $\epsilon$.

The multiplicative constant has an explicit value, that requires some
more notation, and does involve hyperbolic geometry. We denote by
$(\Omega_i)_{i\in I}$ a family of representatives, modulo the action
of $\Ga$, of the connected components of $\Omega\Ga$ whose stabilisers
have infinite index in $\Ga$. For every $i\in I$, let $\C\Omega_i$ be
the convex hull of $\Omega_i$ in the upper-half space model of the
$3$-dimensional real hyperbolic space $\HH^3_\RR$, and let
$\sigma^-_{\C\Omega_i}$ be the (inner) skinning measure of
$\C\Omega_i$ for $\Ga$. We also denote by $\operatorname{HB}_\infty$
the horoball in $\HH^3_\RR$ consisting of points with vertical
coordinates at least $1$, and by $\sigma^+_{\operatorname{HB}_\infty}$
its (outer) skinning measure  for $\Ga$.

\bcoro [Parkkonen-Paulin] \label{coro:compdomdiscont} Let $\Ga$ be a
geometrically finite discrete group of \\$\PSL_{2}(\CC)$, with bounded
and not totally disconnected limit set in $\CC$, which does not
contain a quasifuchsian subgroup with index at most $2$. Assume that
the Hausdorff dimension $\delta$ of the limit set of $\Ga$ is at least
$\frac{1}{2}$. Then there exists $\kappa>0$ such that the number of
connected components of the domain of discontinuity $\Omega\Ga$ of
$\Ga$ with diameter at least $\epsilon$ is equal, as $\epsilon\ra 0$,
to
$$
\frac{2^\delta\;\|\sigma^+_{\operatorname{HB}_\infty}\|\;
\sum_{i\in I}\|\sigma^-_{\C\Omega_i}\|}{\delta\;\|m_{\rm BM}\|}
\,\epsilon^{-\delta} +\bigO(\epsilon^{-\delta+\kappa})\;.
$$  
\ecoro

We refer to \cite[Coro.~25]{ParPau13b} for a proof of this result,
with the error term coming from \cite[Theo.~28]{ParPau13b} and
\cite{MohOh14}, as explained in the discussion of Assertion (4) of
Theorem \ref{theo:requisitgeodflow}. This corollary largely extends
the result of Oh-Shah \cite{OhShaCircles} when all the connected
components of the domain of discontinuity are assumed to be round
discs. Note that the fractal geometry of the boundary of general
convex hulls is an important feature that needs to be adressed by
non-homogeneous dynamics arguments.

%trim option's parameter order: left bottom right top
\begin{center}
\includegraphics[trim = 27mm 90mm 22mm 110mm, clip, width=0.8 \textwidth]{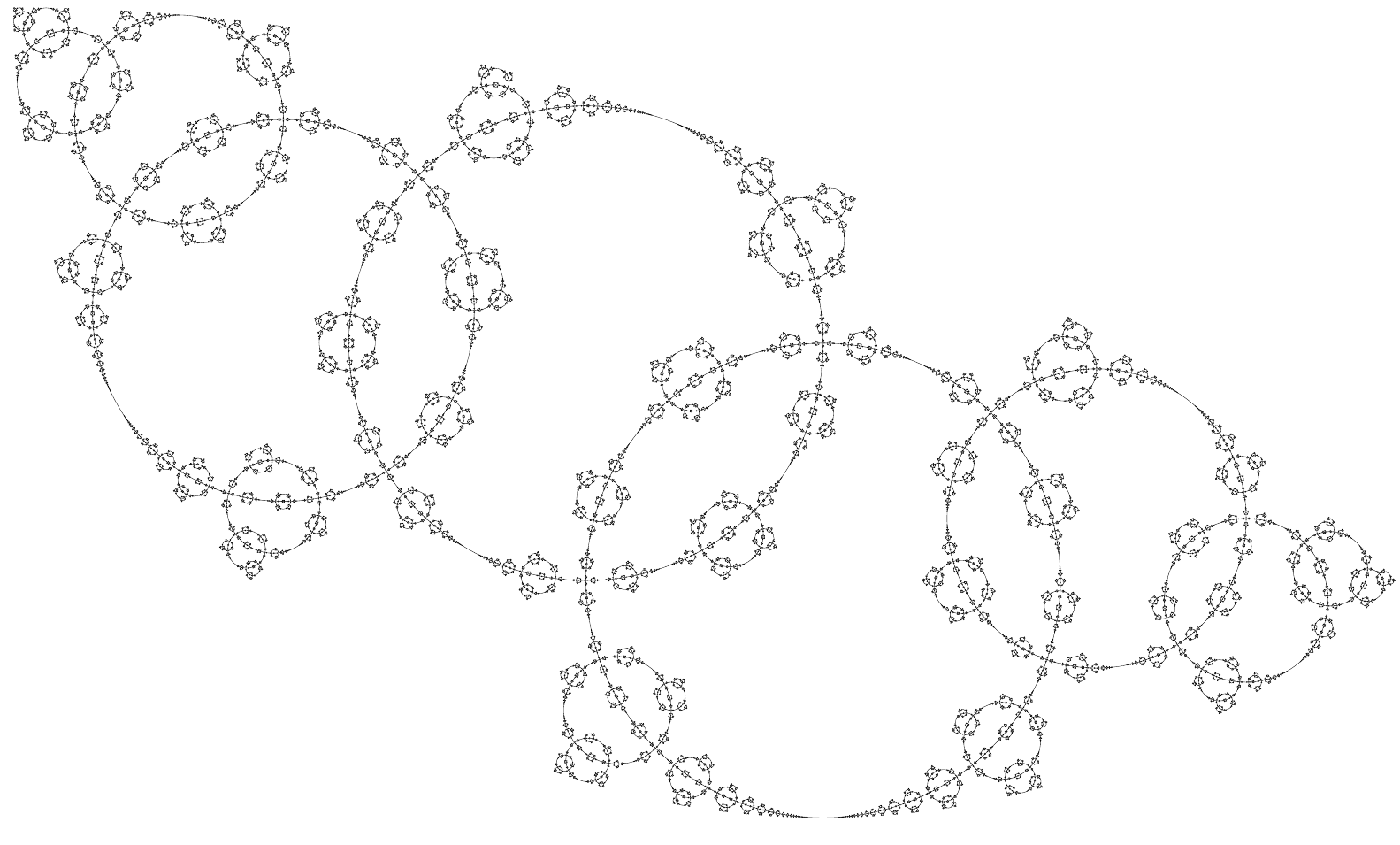}
\end{center}

\subsection{Equidistribution and counting of common perpendicular}
\label{subsect:commonperp}

Let $(\wt M, \Ga, M)$ be as in the beginning of Subsection
\ref{subsec:negacurv}. A {\it common perpendicular} from a closed
convex subset $A^-$ in $\wt M$ to a closed convex subset $A^+$ in $\wt
M$ is a geodesic arc $\wt \alpha$ in $\wt M$ whose initial tangent
vector $v^-_{\wt \alpha}$ belongs to $\normalout A^-$ and terminal
tangent vector $v^+_{\wt \alpha}$ belongs to $\normalin A^+$. Note
that there exists such a common perpendicular if and only if $A^-$ and
$A^+$ are nonempty, proper, with disjoint closures in $\wt
M\cup\partial_\infty\wt M$. It is then unique, and its length is
positive.

Let $C^\pm$ be nonempty proper closed convex subsets of $\wt M$, with
stabiliser $\Ga_{C^\pm}$ in $\Ga$, such that the family $(\ga
C^\pm)_{\ga\in \Ga/\Ga_{C^\pm}}$ of subsets of $\wt M$ is locally
finite. We denote by $\Perp(C^-,C^+)$ the set of images in $M=\Ga\bs
\wt M$ of the common perpendiculars from $\ga^- C^-$ to $\ga^+ C^+$ as
$\ga^\pm$ ranges over $\Ga$, and, for every $t>0$, by
$\Perp(C^-,C^+,t)$ the subset of the ones with length at most $t$. For
every $\alpha\in\Perp(C^-,C^+)$, we denote by $v^-_{\alpha}$ and
$v^+_{\alpha}$ its initial tangent vector and terminal tangent vector,
which belong to the image in $T^1M$ of respectively $\normalout C^-$
and $\normalin C^+$.

\medskip Since $\Ga$ might have torsion, and since $\Ga_{C^\pm}\bs
C^\pm$ does not necessarily embed in $M=\Ga\bs \wt M$, each element
$\alpha$ of $\Perp(C^-,C^+,t)$ comes with a natural multiplicity
$m(\alpha)$.  Denote by $\wt\alpha$ any common perpendicular in $\wt
M$ between translates of $C^-$ and $C^+$ with image $\alpha$ in $M$
and by $\Ga_{\wt\alpha}$ its pointwise stabiliser in $\Ga$. Then
$$
m(\alpha)=\frac{\card\{\ga^-\in\Ga/\Ga_{C^-}\;:\;
v^-_{\wt \alpha}\in\ga^-\normalout C^-\}\;\card\{\ga^+\in\Ga/\Ga_{C^+}\;:\;
v^+_{\wt \alpha}\in\ga^+\normalin C^+\}}{\card \;\Ga_{\wt\alpha}}\;.
$$
Note that the numerator and the denominator are finite by the local
finiteness of the families $(\ga C^\pm)_{\ga\in \Ga/\Ga_{C^\pm}}$ and
the discreteness of $\Ga$, and they depend only on the orbit of $\wt
\alpha$ under $\Ga$. This multiplicity is indeed natural. Concerning
the denominator, in any counting problem of objects possibly having
symmetries, the appropriate counting function consists in taking as
the multiplicity of an object the inverse of the cardinality of its
symmetry group. The numerator is here in order to take into account
the fact that the elements $\ga^\pm$ in $\Ga$ such that $\wt \alpha$
is a common perpendicular from $\ga^- C^-$ to $\ga^+ C^+$ are not
necessarily unique, even modulo $\Ga_{C^\pm}$.  The natural counting
function of the common perpendiculars between the images of $C^\pm$ is
then the map
$$
t\mapsto \N_{C^-,\,C^+}(t) =
\sum_{\alpha\in\Perp(C^-,\,C^+,\,t)}\;m(\alpha)\;.
$$
The reader can assume for simplicity that $\Ga$ is torsionfree and
that $\Ga_{C^\pm}\bs C^\pm$ embeds in $M=\Ga\bs \wt M$ by the map
induced by the inclusion of $C^\pm$ in $M$, in which case all
multiplicities are $1$.

\medskip Below, we state our equidistribution result, in the outer and
inner unit normal bundles of the images in $M$ of $C^-$ and $C^+$, of
the initial and terminal tangent vectors of the common perpendiculars
between the images of $C^-$ and $C^+$, as well as our asymptotic
formula as $t\ra+\infty$ for the number of common perpendiculars of
length at most $t$ between the images of $C^-$ and $C^+$. We refer
respectively to \cite[Theo.~14, 28]{ParPau13b} and \cite[Coro.~20,
28]{ParPau13b} for more general versions, involving more general
locally finite families of convex subsets, versions with weights
coming from potentials (the Bowen-Margulis measure being then replaced
by the Gibbs measure of \cite{PauPolSha}), and for version with error
terms under an additional assumption of exponential decay of
correlations (see Theorem \ref{theo:requisitgeodflow} (4)).

\btheo[Parkkonen-Paulin] \label{theo:mainequidis} Assume that $m_{\rm
  BM}$ is finite and mixing under the geodesic flow, and that the
skinning measures $\sigma^\pm_{C^\mp}$ are finite.  For the weak-star
convergence of measures on $T^1M\times T^1M$, we have
$$
\lim_{t\ra+\infty}\; \delta_\Ga\;\|m_{\rm BM}\|\;e^{-\delta_\Ga t}
\sum_{\alpha\in\Perp(C^-,\,C^+,\,t)} \;m(\alpha)\; 
\Delta_{v^-_\alpha} \otimes\Delta_{v^+_\alpha}\;=\;
\sigma^+_{C^-}\otimes \sigma^-_{C^+}\,.
$$
\etheo

\btheo[Parkkonen-Paulin] \label{theo:maincount} Assume that $m_{\rm
  BM}$ is finite and mixing under the geodesic flow, and that the
skinning measures $\sigma^\pm_{C^\mp}$ are finite and nonzero. Then,
as $t\to+\infty$,
$$
\N_{C^-,\,C^+}(t)\sim
\frac{\|\sigma^+_{C^-}\|\,\|\sigma^{-}_{C^+}\|}
{\|m_{\rm BM}\|}\,\frac{e^{\delta_{\Ga} \,t}}{\delta_{\Ga}}\,.
$$
\etheo

The counting function $\N_{C^-,\,C^+}(t)$ has been studied in
various special cases since the 1950's and in a number of recent
works, sometimes in a different guise, see the survey \cite{ParPauRev}
for more details. A number of special cases were known before our
result:
\begin{itemize}
\item $C^-$ and $C^+$ are reduced to points, by for instance
  \cite{Huber59}, \cite{Margulis69} and \cite{Roblin03},
\item $C^-$ and $C^+$ are horoballs, by \cite{BelHerPau01},
  \cite{HerPau04}, \cite{Cosentino99} and \cite{Roblin03} without an
  explicit form of the constant in the asymptotic expression,
\item $C^-$ is a point and $C^+$ is a totally geodesic submanifold, by
  \cite{Herrmann62}, \cite{EskMcMul93} and \cite{OhShaCircles} in
  constant curvature,
\item $C^-$ is a point and $C^+$ is a horoball, by
  \cite{Kontorovich09} and \cite{KonOh11} in constant curvature, and
  \cite{Kim13} in rank one symmetric spaces,
\item $C^-$ is a horoball and $C^+$ is a totally geodesic submanifold,
  by \cite{OhSha12} and \cite{ParPau12JMD} in constant curvature, and
\item $C^-$ and $C^+$ are (properly immersed) locally geodesic lines
  in constant curvature and dimension $3$, by \cite{Pollicott11}.
\end{itemize}

As a new particular case, if $M$ has constant curvature $-1$, if the
images in $M$ of $C^-$ and $C^+$ are closed geodesics of lengths
$\ell_-$ and $\ell_+$, respectively, then the number of common
perpendiculars (counted with multiplicity) from the
image of $C^-$ to the
image of $C^+$ of
length at most $s$ satisfies, as $s\ra+\infty$,
\begin{equation}\label{eqcaseclosedgeod}
\N_{C^-,\,C^+}(s)\sim
\frac{\pi^{\frac{n}{2}-1}\Ga(\frac{n-1}{2})^2}
{2^{n-2}(n-1)\Ga(\frac{n}{2})}
\;\frac{\ell_-\ell_+}{\Vol(M)}\;
e^{(n-1)s}\;.
\end{equation}
When $M$ is a closed hyperbolic surface and $C^-=C^+$, this formula
\eqref{eqcaseclosedgeod} has been obtained by Martin-McKee-Wambach
\cite{MarMcKWam11} by trace formula methods, though obtaining the case
$C^-\neq C^+$ seems difficult by these methods.

\medskip The family $(\ell(\alpha))_{\alpha\in\Perp(C^-,\,C^+)}$ (with
multiplicities) will be called the {\it marked ortho\-length spectrum}
from $C^-$ to $C^+$. The set of lengths (with multiplicities) of
elements of $\Perp(C^-,C^+)$ will be called the {\it ortholength
  spectrum} of $C^-,C^+$. This second set has been introduced by
Basmajian \cite{Basmajian93} (under the name ``full orthogonal
spectrum'') when $M$ has constant curvature, and the images in $M$ of
$C^-$ and $C^+$ are disjoint or equal embedded totally geodesic
hypersurfaces or embedded horospherical cusp neighbourhoods or
embedded balls (see also \cite{BriKah10} when $M$ is a compact
hyperbolic manifold with totally geodesic boundary and the images in
$M$ of $C^-$ and $C^+$ are exactly $\partial M$). The two results
above are hence major contributions to the asymptotics of marked
ortholength spectra.

\medskip 
Let us give a brief sketch of the proof of Theorem
\ref{theo:maincount}, refering to \cite[\S 4.1]{ParPau13b} for a full
proof.

\medskip\noindent{\bf Step 1. } In this technical step, we start by
constructing dynamical neighbourhoods and test functions around the
outer/inner unit normal bundles of our convex sets, that will be
appropriately pushed forward/backwards by the geodesic flow (using
the nice contraction/dilation properties of Hamenst\"adt's distances
compared to the ones of the induced Riemannian metric in variable
curvature). We fix $R>0$ big enough, and we will let $\eta>0$, a
priori small enough, tend to $0$.

\medskip\noindent
\begin{minipage}{8.9cm}
  ~~~ For all $w\in T^1\wt M$, let $V^+_{w}$ be the open ball of
  center $w$ and radius $R$ for Hamenst\"adt's distance on the strong
  stable leaf $W^+(w)$ of $w$. For every $\eta>0$, let
$$
\V^+(C^-)=
\bigcup_{w\in \normalout C^-,\; s\in\,]-\eta,\eta[} \flow s V^+_{w}\,,
$$
which is a neighbourhood of $\normalout C^-$.
\end{minipage}
\begin{minipage}{6cm}
\begin{center}
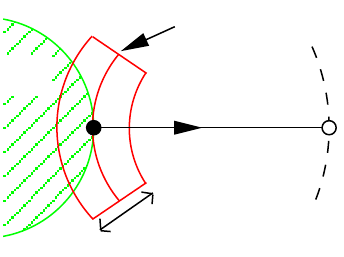
\end{center}
\end{minipage}

Let $h^-: T^1\wt M\ra \mathopen{]}0,+\infty]$ be the measurable and
$\wt m_{\rm BM}$-almost everywhere finite map (since its denominator
is positive if $w^-\in\Lambda\Ga$) defined by
$$
h^-(w)=\frac{1}{2\eta\;\wt\sigma^-_{\operatorname{HB}^+(w)}(V^+_{w})}\;.
$$

Let us denote by $\mathbbm{1}_A$ the characteristic function of a
subset $A$. Consider the map $\wt\psi_\eta^-: T^1\wt M\to[0,+\infty]$
defined by
\begin{equation}\label{eq:defiphi}
\wt\psi_\eta^-:v\mapsto\sum_{\ga\in\Ga/\Ga_{C^-}}\;h^-\circ f^+_{\ga C^-}(v)\;\;
\mathbbm{1}_{\V^+(\ga C^-)}(v)\,.
\end{equation}
This map is $\Ga$-invariant and measurable, hence it defines a measurable
map $\psi_\eta^-:T^1M\to[0,+\infty]$, with support in a neighbourhood of
the image of $\partial^1_+C^-$ in $T^1M$. We define similarly
$\psi_\eta^+:T^1M\to[0,+\infty]$ with support in a neighbourhood of the
image of $\partial^1_-C^+$ in $T^1M$.

\medskip By the disintegration result of Proposition
\ref{prop:disintegration}, the functions $\psi_\eta^\mp$ are
integrable and
\begin{equation}\label{eq:corodisint}
  \int_{T^1 M} \psi_\eta^\mp\;dm_{\rm BM}=\|\sigma^\pm_{C^\mp}\|\,.
\end{equation}

\medskip\noindent{\bf Step 2. } In this step, we use the mixing
property of the geodesic flow, as first introduced by Margulis in his
thesis (see for instance \cite{Margulis04a}). Due to the symmetry of
the problem, a one-sided pushing of the geodesic flow, as in all the
previous works using Margulis's ideas, is not sufficient, and we need
to push simultaneously the outer and inner unit normal vectors to the
convex sets in opposite directions.

For all $t\geq 0$, let
$$
a_\eta(t)= \int_{T^1 M} \psi_\eta^-\circ \flow{-t/2}\;\;
\psi_\eta^+\circ\flow{t/2}\;dm_{\rm BM}\,.
$$
Then the mixing hypothesis of the geodesic flow and Equation
\eqref{eq:corodisint} imply that
$$
\lim_{t\ra+\infty} a_\eta(t)=
\frac{\|\sigma^+_{C^-}\|\;\|\sigma^-_{C^+}\|}{\|m_{\rm BM}\|}\,.
$$

\medskip\noindent{\bf Step 3. } In this step, we give another estimate
of $a_\eta(t)$, relating it to the counting of common
perpendiculars. One of the main new ideas in the proof (see \cite[\S
2.3]{ParPau13b} for a complete version) is an effective study of the
geometry and the dynamics of the accidents that occur around midway of
the pushing by the geodesic flow, yielding an effective statement of
creation of common perpendiculars.

\begin{center}
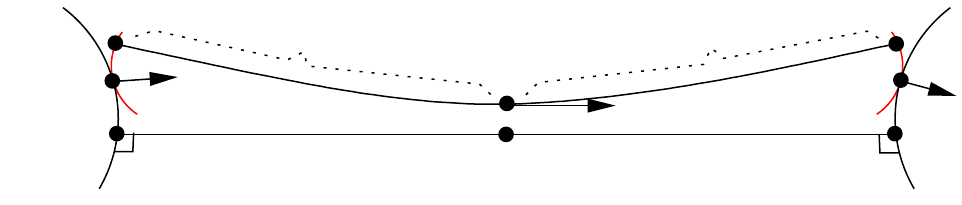
\end{center}

In order to give an idea of this phenomenon, assume that $v\in T^1M$
belongs to the support of the function $\psi_\eta^-\circ \flow{-t/2}\;\;
\psi_\eta^+ \circ\flow{t/2}$ whose integral is $a_\eta(t)$. This is
equivalent to assuming that $\flow{\pm t/2}v$ belongs to the support
of $\psi_\eta^\pm$. If $\wt v$ is a lift of $v$ to $T^1\wt M$, by the
definition of the maps $\psi_\eta^\mp$, this is equivalent to asking that
there exist $\ga^\pm\in\Ga$, $w^\pm\in\normalmp C^\pm$ and $s^\pm\in
\mathopen{]}-\eta, \eta\mathclose{[}$ such that $\flow{\pm
  (t/2+s^\pm)}\wt v$ belongs to Hamenstädt's balls $V^\mp_{w^\pm}$.
If $R$ is fixed, when $\eta>0$ is small enough and $t$ is big enough,
negative curvature estimates say that $v$ is very close to the tangent
vector at the midpoint of a common perpendicular from $\ga^-C^-$ to
$\ga^+C^+$ whose length is close to $t$ (see \cite[\S
2.3]{ParPau13b}).

To obtain a precise estimate on $a_\eta(t)$, given a fundamental
domain $\F$ for the action of $\Ga$ on $T^1\wt M$, we apply Fubini's
theorem, as in Sarnak's ``unfolding technique'':
$$
\int_{\F}\sum_{\ga^-\in\Ga/\Ga_{C^-}}\sum_{\ga^+\in\Ga/\Ga_{C^+}}=
\sum_{\ga^-\in\Ga/\Ga_{C^-}}\sum_{\ga^+\in\Ga/\Ga_{C^+}}\int_{\F}\;,
$$
as well as a fine analysis (especially refined for the error term
estimates) of the intrinsic geometry in variable curvature (almost
everywhere defined second fundamental form, ...)  of the outer/inner
unit normal bundle pushed a long time by the geodesic flow.  We then
conclude by a Cesaro-type of argument in order to consider all common
perpendiculars with length at most $T$, as $T$ tends to $+\infty$, and
by letting $\eta$ tend to $0$.

\subsection{Towards the arithmetic applications}
\label{subsect:fromgeomtoarith}

As we already hinted to, Theorems \ref{theo:equidistrace},
\ref{theo:relatcomplcount}, \ref{theo:equidistribR},
\ref{theo:equidistribC}, \ref{theo:equidisHeis}, \ref{theo:countchain}
and \ref{theo:equidischain} all follow from Theorem
\ref{theo:mainequidis} or Theorem \ref{theo:maincount}, though many
more tools and ideas are needed, in particular volume computations of
arithmetic orbifolds.

We only indicate the very beginning of the proof of these theorems,
giving a bit more details on Theorems \ref{theo:Lagrangequadirrat} to
\ref{theo:equidistribR}.

\medskip To prove Theorem \ref{theo:equidistribC}, we apply Theorem
\ref{theo:mainequidis} with $\wt M=\HH^3_\RR$, $\Ga$ the Bianchi group
$\PSL_2(\OOO_K)$, and $C^-=C^+$ any horoball centered at $\infty$ in
the upper halfspace model of $\HH^3_\RR$. Note that the cusps of the
noncompact finite volume hyperbolic orbifold $\PSL_2(\OOO_K)\bs
\HH^3_\RR$ correspond to the ideal classes of $\OOO_K$ (in particular
if $K=\QQ(i)$, there is only one cusp), see for instance
\cite{ElsGruMen98}. Keeping the same $C^-$ and taking $C^+$ centered
at a parabolic fixed point defining the cusp allows version of Theorem
\ref{theo:equidistribC} when $p$ and $q$ are varying in a given
fractional ideal of $\OOO_K$ (when the class number of the imaginary
quadratic number field $K$ is larger than $1$).

\medskip To prove Theorem \ref{theo:equidisHeis}, we consider the
Hermitian form $h:(z_0,z_1,z_2)\mapsto -z_0\overline{z_2}-
z_2\overline{z_0} +|z_1|^2$ on $\CC^3$ whose signature is $(1,2)$. We
apply Theorem \ref{theo:mainequidis} with $\wt M$ the projective model
$\{[z_0:z_1:z_2]\in \PP_2(\CC)\;:\; h(z_0,z_1,z_2)<0\}$ of the complex
hyperbolic plane $\HH^2_\CC$, with $\Ga$ the Picard group
$\PSU_h(\OOO_K) = \PSU_h\cap \operatorname{PGL}_3(\OOO_K)$, and with
$C^-=C^+$ any horoball centered at $\infty=[1:0:0]$. See
\cite[Theo.~12, 13]{ParPauHeis} for the other ingredients of the
proof. The reason why large white regions appear around the points
$\frac{p}{q}\in K$ with $|q|$ small in the figures before and after
Theorem \ref{theo:equidistribC} is that the horoballs in the
$\Ga$-orbit of $C^-$ centered at these points have large Euclidean
radius, hence it is (quadratically) difficult to fit disjoint
horoballs in this orbit below them.

Replacing in the above data $C^+$ by the convex hull in $\HH^2_\CC$ of
an arithmetic chain $C_0$, applying Theorem \ref{theo:maincount} and
\ref{theo:mainequidis} is the very first step for proving Theorem
\ref{theo:countchain} and \ref{theo:equidischain}, respectively. See
\cite[Theo.~19, 20]{ParPauHeis} for the other ingredients of the
proof.

\medskip To prove Theorems \ref{theo:equidistrace},
\ref{theo:relatcomplcount} and \ref{theo:equidistribR}, we apply
Theorem \ref{theo:mainequidis} or Theorem \ref{theo:maincount} with
$\wt M$ the upper halfplane model of $\HH^2_\RR$ and $\Ga=\PSL_2(\ZZ)$
(or appropriate finite index subgroups when we want versions with
additional congruence assumptions). Note that the modular curve
$\PSL_2(\ZZ)\bs \HH^2_\RR$, being arithmetic hyperbolic, has its
Bowen-Margulis measure proportional to its Liouville measure, hence
finite, and its geodesic flow is mixing, with exponential decay of
correlation.  We then take
\begin{itemize}
\item $C^-$ a horoball centered at $\infty$ and $C^+$ the
geodesic line $]\alpha_0,\alpha_0^\sigma[$ in $\HH^2_\RR$ with points
at infinity $\alpha_0,\alpha_0^\sigma$ for Theorem
\ref{theo:equidistrace},

\item $C^-$ the geodesic line $]\alpha_0,\alpha_0^\sigma[$ and
$C^+$ the geodesic line $]\beta_0,\beta_0^\sigma[$ for
Theorem \ref{theo:relatcomplcount},

\item $C^-$ the geodesic line $]\alpha_0,\alpha_0^\sigma[$ and
$C^+$ a horoball centered at $\infty$ for Theorem 
\ref{theo:equidistribR}.
\end{itemize}

\begin{center}
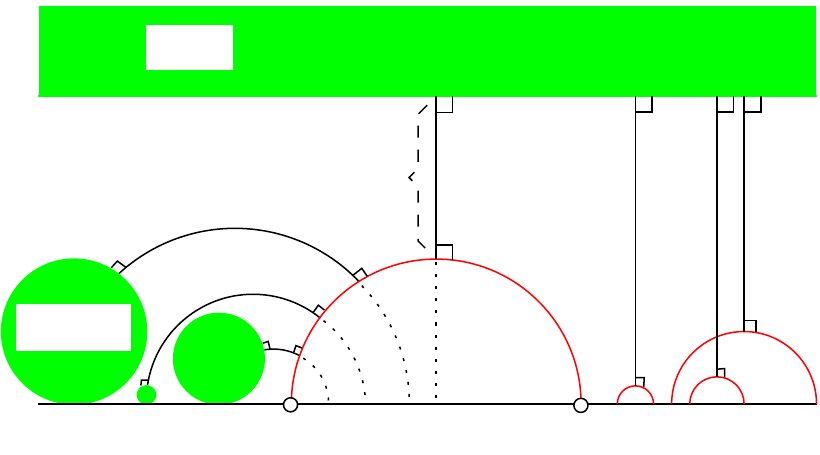
\end{center}

The key input, also crucial to prove Theorem
\ref{theo:Lagrangequadirrat}, is the well-known hyperbolic geometry
understanding of quadratic irrationals. A real number $\alpha$ is a
quadratic irrational if and only if it is fixed by (the action by
homography of) an element $\ga\in\PSL_2(\ZZ)$ with $|\operatorname{tr}
\; \ga|>2$. Then $\alpha^\sigma$ is the other fixed point of $\ga$,
and the geodesic line $L_\alpha= \mathopen{]}\alpha, \alpha^\sigma
\mathclose{[}$ maps to a closed geodesic in the hyperbolic orbifold
$\PSL_2(\ZZ)\bs\HH^2_\RR$. In particular, the image of $\normalpm
L_\alpha$ in $\PSL_2(\ZZ) \bs T^1\HH^2_\RR$ is compact, and the
skinning measures $\sigma^\pm_{L_\alpha}$ are positive and finite.

The first hint that there is a connection between quadratic
irrationals and common perpendiculars is the following one. Let
$\operatorname{HB}_\infty$ be the horoball centered at $\infty$,
consisting of the points of $\HH^2_\RR$ with Euclidean height at least
$1$. Note that its stabiliser in $\PSL_2(\ZZ)$ acts cocompactly on
$\normalpm\operatorname{HB}_\infty$, hence the skinning measures
$\sigma^\pm_{\operatorname{HB}_\infty}$ are positive and finite.
Then, by an easy computation in hyperbolic geometry, the common
perpendicular between $\operatorname{HB}_\infty$ and the geodesic line
$]\alpha,\alpha^\sigma[$ (assuming that they are disjoint) has length
$\log H(\alpha)$, where
$$
H(\alpha)=\frac{2}{|\alpha-\alpha^\sigma|}\;.
$$

Another important observation to prove Theorem \ref{theo:equidistribR}
(taking $\alpha=\frac{1+\sqrt{5}}{2}$ the Golden Ratio,
$C^-=L_{\alpha}$ and $C^+=\operatorname{HB}_\infty$) is that since
the modular curve has finite volume, the skinning measure on
$\normalout C^-$ is homogeneous. Hence, on each of the two connected
components of $\normalout C^-$ which are naturally parametrised by
$\RR$, it is proportional to the Lebesgue measure, and this Lebesgue
measure projects by the negative endpoint map $v\mapsto v^-$ to a
measure on $\RR-\{\alpha,\alpha^\sigma\}$ proportional to
$\frac{d\Leb_\RR(t)}{|t^2-t-1|}$, explaining the limit in
Theorem \ref{theo:equidistribR}. See \cite[Theo.~6]{ParPau14AFST} for a
complete proof.

{\small\bibliography{../biblio} }
%{\small \bibliography{../viitteet} }

\bigskip
{\small
\noindent \begin{tabular}{l} 
Department of Mathematics and Statistics, P.O. Box 35\\ 
40014 University of Jyv\"askyl\"a, FINLAND.\\
{\it e-mail: jouni.t.parkkonen@jyu.fi}
\end{tabular}
\medskip

\noindent \begin{tabular}{l}
D\'epartement de math\'ematique, UMR 8628 CNRS, B\^at.~425\\
Universit\'e Paris-Sud,
91405 ORSAY Cedex, FRANCE\\
{\it e-mail: frederic.paulin@math.u-psud.fr}
\end{tabular}
}

\end{document}